\documentclass[review]{elsarticle}

\usepackage{lineno,hyperref}
\modulolinenumbers[5]

\journal{Journal of Computational Physics}

\usepackage{graphicx}
\usepackage{amssymb}
\usepackage{amsmath}
\usepackage{epstopdf}
\usepackage{url}

\newcommand\Beq{\begin{eqnarray}} 
\newcommand\Eeq{\end{eqnarray}}
\newcommand{\n}{  \\ }

\newcommand{\eq}[1]{equation~(\ref{#1})}
\newcommand{\eqs}[2]{equations~(\ref{#1})~\&~(\ref{#2})}
\newcommand{\eqss}[2]{equations~(\ref{#1})--(\ref{#2})}
\newcommand{\Eq}[1]{Equation~(\ref{#1})}
\newcommand{\Eqs}[2]{Equations~(\ref{#1})~\&~(\ref{#2})}

\newcommand{\Fig}[2]{Figure~\ref{#1}{\textit{#2}}}

\newcommand{\ie}{\textit{i.e.}, }
\newcommand{\eg}{\textit{e.g.}, }

\newcommand{\dd}[1]{\,\mathrm{d}{#1}}
\newcommand{\bra}[1]{\big<\, #1 \, \big|}
\newcommand{\ket}[1]{\big|\, #1 \, \big>}
\newcommand{\braket}[2]{\big<\, #1 \, \big| \, #2 \, \big>}

\begin{document}

\begin{frontmatter}

\title{Tensor calculus in polar coordinates using Jacobi polynomials}

\author[sydney]{Geoffrey M. Vasil\corref{mycorrespondingauthor}}
\cortext[mycorrespondingauthor]{\textit{Corresponding author email:} \url{geoffrey.vasil@sydney.edu.au}}

\author[MIT]{Keaton J.~Burns}

\author[berkeley]{Daniel Lecoanet}

\author[sydney]{Sheehan Olver}

\author[CU]{Benjamin P.~Brown}

\author[bates]{Jeffrey S.~Oishi}

\address[sydney]{School of Mathematics \& Statistics, University of Sydney, NSW 2006, Australia}
\address[MIT]{Dept.~Physics, Massachusetts Institute of Technology, Cambridge, MA 02139, USA}
\address[berkeley]{Dept.~Physics and Theoretical Astrophysics Center, University of California, Berkeley, CA 94720, USA}
\address[CU]{LASP and Dept.~Astrophysical \& Planetary Sciences, University of Colorado, Boulder, CO 80309, USA}
\address[bates]{Dept.~Physics \& Astronomy, Bates College, Lewiston, ME 04240, USA}

\begin{abstract}
Spectral methods are an efficient way to solve partial differential equations on domains possessing certain symmetries.  The utility of a method depends strongly on the choice of spectral basis.  In this paper we describe a set of bases built out of Jacobi polynomials, and associated operators for solving scalar, vector, and tensor partial differential equations in polar coordinates on a unit disk.  By construction, the bases satisfy regularity conditions at $r=0$ for any tensorial field. The coordinate singularity in a disk is a prototypical case for many coordinate singularities. The work presented here extends to other geometries. The operators represent covariant derivatives, multiplication by azimuthally symmetric functions, and the tensorial relationship between fields.  These arise naturally from relations between classical orthogonal polynomials, and form a Heisenberg algebra. Other past work uses more specific polynomial bases for solving equations in polar coordinates. The main innovation in this paper is to use a larger set of possible bases to achieve maximum bandedness of linear operations. We provide a series of applications of the methods, illustrating their ease-of-use and accuracy.
\end{abstract}


\end{frontmatter}

\section{Introduction}

Cylindrical polar coordinates find applications in countless areas of science and engineering. Important applications include pipe flow, laboratory studies of thermal convection, astrophysical accretion disks, electromagnetic waveguides, elastic deformation of rods, astronomical instrumentation, and plasma tokamaks.  Many applications require the accurate and efficient solution of systems of partial differential equations (PDEs).  Pseudospectral methods of different types prove useful for this task in many different geometries. In polar coordinates, the periodic nature of the azimuth angle allows the effective use of Fourier series, where
\Beq
\label{Fourier series}
f(r,\theta) = \sum_{m=-\infty}^\infty f_m(r) e^{i m \theta},\quad 
f_{m}(r)  = \frac{1}{2\pi}\int_{0}^{2\pi} f(r,\theta)e^{-im \theta} \dd{\theta}.
\Eeq
After the Fourier transform, differentiation in $\theta$ becomes multiplication by $i m$.  While Fourier analysis easily dispatches the azimuthal coordinate for functions on a disk, the radial coordinate presents difficulty for the following reason. For functions analytic everywhere on the disk, including the origin,   
\Beq
\label{analytic}
f_{m}(r) \sim r^{m} F(r^2) \quad \mathrm{as} \quad r \to 0,
\Eeq
where $F(r^2)$ is an even function of $r$ that is analytic at the origin. 

The coordinate singularity at the disk centre requires an $m$-th order zero for infinite differentiability \cite{orszag_1974,boyd2011comparing}.
Enforcing this condition in numerical calculations presents challenges; especially for $m \gg 1$.  Many authors address this challenge with equally as many different techniques. 
Even considering regularity at the origin, the disk geometry allows a large number of possible 
orthogonal-polynomial bases \cite{dunklxu,koornwinder}.  Zernike (1934) \cite{zernike}  produced the first practical set of polynomials for expanding functions on the unit disk. This basis proves particularly useful in optical applications.  Bhatia \& Wolf (1954) \cite{bhatia1954circle} pointed out that this set is the only out of a possible infinity that contains \textsl{``simple properties strictly analogous to that of Legendre polynomials.''}

Boyd \& Yu (2011) \cite{boyd2011comparing} provide a comprehensive review of the history and contemporary methods used to solve Poisson's equation in a disk. In particular, the paper reviews bases using Zernike-type polynomials, as well as the more common Chebyshev polynomials.  The results for Chebyshev series range from acceptable to untenable. The diversity of  Chebyshev methods results from different ways to represent the pole condition and/or the reflectional symmetry near the origin. A minimalist approach happens to produce the best option. This option expands even/odd-$m$ modes in terms of an even/odd-degree Chebyshev series. This approach double wraps the disk using a Chebyshev series over $-1 \le r \le 1$ \cite{trefethenspectralmethods,fornberg1998practical}. Compared to other Chebyshev options, simple even-odd matching works well with no other special intervention \cite{livermore_jones_worland_2007}. Even-odd matching and/or double-covering can satisfy \eq{analytic} with good-to-moderate accuracy. These schemes however do not enforce the analytic condition explicitly. This implies that singularities can still arise in higher-order derivatives; also see \cite{livermore_jones_worland_2007}.  Even weak singularities can produce instabilities at the origin when performing time-evolution simulations.  As a third option, the Roberts basis combines an even Chebyshev series with an explicit $r^{m}$ prefactor. In spite of initial attractiveness (\eg possessing a fast transform), this basis suffers from extreme numerical ill conditioning, and is not recommended \cite{boyd2011comparing}.

Regarding the Zernike-type bases, Boyd \& Yu point out that they are \textsl{``More accurate for large m.''} They also discuss the less-fortunate fact that Zernike bases do not admit a fast transform in the radial direction. But that for various reasons \textsl{``the advantages of `FFT-ability' is not huge.''}. 
They conclude that \textsl{``It is difficult to definitely endorse one particular method for the disk because of the vast diversity of solutions to interesting engineering and science problems.''} Furthermore, Slevinsky (2016) has recently made significant progress toward designing an effective fast transform from values at Chebyshev points to Jacobi coefficients that would work with the Zernike basis \cite{slevinsky_2016}. For these reasons and more, we believe that polynomial bases that satisfy \eq{analytic} are very useful in many applications and are worthy of more detailed understanding.   

In addition to scalar-valued functions, many situations also require vector and tensor fields. Vectors introduce additional complications near the coordinate singularity. Much less work exists addressing these issues.  In particular, the $m$-th Fourier components of a vector field behave such that
\Beq
\label{v-analytic}
v_{m}(r) \sim r^{m-1}V(r^2) \quad r \to 0,
\Eeq
where, like $F(r^2)$ in \eq{analytic}, $V(r^2)$ is an even function of $r$ that is analytic at the origin. We can (for example) see the necessity of \eq{v-analytic} by differentiating \eq{analytic} with respect to $r$.

Sakai \& Redekopp (2009) \cite{sakai2009application} circumvented this issue by working with rescaled variables of the form $r v_{m}(r)$; which behaves like \eq{analytic}.  Li, Livermore \& Jackson (2010) \cite{li2010optimal} use a poloidal-toroidal formulation to create a genuine (higher-order) scalar system out of a specific vector system.  Using a technique equivalent to the $r$ rescaling, Matsushima \& Marcus (1995) \cite{matsushima1995spectral} show (and Boyd \& Yu \cite{boyd2011comparing} reiterate) that Zernike polynomials produce \textsl{pentadiagonal} matrices for the solution of the radial portion of Poisson's equation. Lastly, Townsend, Wilber \& Wright (2016) develop an efficient low-rank approximation of scalar and vector functions on the disk that preserve regularity \cite{townsend_etal_2016}. These methods work well for data analysis. Their application to time evolving systems remains less clear. 
 
We show in this paper the non-necessity of radial rescaling and/or equation reformulation. Previous works found recursion relationships for elementary operators $r^2$ and $r d/dr$. Our calculus finds simpler factorisations of these operations in terms of $r$, $d/dr$ and $m/r$. These are all of the elementary operators needed for full tensor calculus. This not only makes calculations easier to formulate, but also more numerically efficient and stable.  In the process, we also show how to construct solutions to Poisson's equation on the disk using only \textsl{tridiagonal} (as opposed to pentadiagonal) matrices; see the discussion in example 1 in \S\ref{examples} for more details.  The foundation of the new results rests on exploiting a more general class of orthogonal polynomials.  That is, we choose different bases to represent domain and range spaces of operators, so that coupling becomes banded.  This mirrors using ultra-spherical polynomials for solving equations on the unit interval \cite{doha2002efficient,doha2009efficient,olver2013fast}. Moreover, we incorporate azimuthally symmetric variable coefficients without destroying bandedness. This occurs via approximating non-constant coefficients with finite-degree polynomials similar to \cite{olver2013fast}. 
 
A central theme of this paper demonstrates that increasing the collection of available bases can increase (\textit{i}) the simplicity of a calculation's numerical implementation; (\textit{ii}) the speed to compute a solution; and (\textit{iii}) the accuracy of the result.  We outline our following results: \S2  derives properties for useful bases for polar coordinates (using properties of Jacobi polynomials). \S3 shows how these bases respond to the covariant derivative operator in polar coordinates. \S4 discusses multiplication by radial functions. \S5 shows how the different bases relate to each other, and how the different operators form a Heisenberg Lie algebra. \S2--5 build the necessary tools to represent and manipulate scalars, vectors, and tensors in different calculations. \S6 applies these tools to a series of four example problems. \S7 gives concluding remarks.  

\section{Calculus of Jacobi polynomials}

Throughout this paper, we use the definition 
\Beq
\label{z-def}
z \equiv 2 r^{2}-1, 
\Eeq
where $0 \le r \le 1$, implies $-1 \le z \le +1$. The two following formulae provide the foundation for a simple calculus of functions on the unit disk.  
For \textsl{any} differentiable function $F(z)$, 
\Beq
\label{m-up}
\left[\frac{d}{dr} - \frac{m}{r} \right] r^{m} F(z) &=& 4r^{m+1} \frac{d}{dz} F(z),\n
\label{m-down}
\left[\frac{d}{dr} + \frac{m}{r} \right] r^{m} F(z) &=& 2r^{m-1} \left[ m  + (1+z) \frac{d}{dz}\right] F(z).
\Eeq
The $z$-differential operators on the right-hand side of \eqs{m-up}{m-down} act on Jacobi polynomials in a particularly simple way. 

Jacobi polynomials $P^{(a,b)}_n(z)$ are the two-parameter set of classical orthogonal polynomials on $-1\le z \le 1$ under the weight $(1-z)^a(1+z)^b$, where $-1 < a,b$:
\Beq
\int_{-1}^{1} P^{(a,b)}_{n}(z)P^{(a,b)}_{n'}(z)\,(1-z)^a\,(1+z)^b \dd{z} \ = \ W^{a,b}_{n} \, \delta_{n,n'},
\Eeq 
and
\Beq
W^{a,b}_{n} \equiv \frac{2^{a+b+1}}{2n+a+b+1} \, \frac{\Gamma(n+a+1) \Gamma(n+b+1)}{\Gamma(n+a+b+1)\,n!},
\Eeq
where $\delta_{n,n'}$ represents the Kronecker delta.
Jacobi polynomials satisfy a Sturm--Liouville differential equation, and a three-term recursion relation. They also satisfy many other useful formulae\footnote{All relevant formulae are found in Digital Library of Mathematical Functions, Chapter on Orthogonal Polynomials, \S18.3--\S18.9: \url{http://dlmf.nist.gov/18}} \cite{DLMF}.  For the purposes of this paper we use: 

For differential operators,
\Beq
2\frac{d}{dz} P_{n}^{(a,b)}(z) &=&  (n+a+b+1) P_{n-1}^{(a+1,b+1)}(z) \label{a+1,b+1},  \\ 
\left[ b  + (1+z) \frac{d}{dz}\right] P_{n}^{(a,b)}(z)  &=&  (n+b)P_{n}^{(a+1,b-1)}(z). \label{a+1,b-1}
\Eeq
For algebraic operators,
\Beq
(1+z)P_{n}^{(a,b)}(z) &=& \tfrac{2(n+1)}{2n + a+ b +1}P_{n+1}^{(a,b-1)}(z) + \tfrac{2(n+b)}{2n + a+ b +1}P_{n}^{(a,b-1)}(z), \label{a,b-1}  \\
P_{n}^{(a,b)}(z) &=& \tfrac{n+a+b+1}{2n + a+ b +1}P_{n}^{(a+1,b)}(z) - \tfrac{n+b}{2n + a+ b +1}P_{n-1}^{(a+1,b)}(z) \label{a+1,b}.
\Eeq
For restriction operators, 
\Beq
\label{P at z=1}
P_{n}^{(a,b)}(z=1) =  \binom{a+n}{n}. \label{z=1} 
\Eeq
Exchanging the roles of $a$ and $b$, we can produce similar relations as \eqss{a+1,b+1}{z=1}, using 
\Beq
\label{exchange}
P_{n}^{(a,b)}(-z) \ = \ (-1)^n P_{n}^{(b,a)}(z).
\Eeq
Doha \& Bhrawy (2006) advocate for using Jacobi polynomials of varying degree to solve a number of problems in one and two dimensions. They specifically state that \textsl{``The key to the efficiency of these algorithms is to construct appropriate base functions, which lead to systems with specially structured matrices that can be efficiently inverted.''}\cite{doha2006efficient}. We expand on this philosophy and construct a hierarchy of basis functions that conform particularly well to vector calculus in circular geometry. 

We define the following radial basis elements in terms of Jacobi polynomials,
\Beq
\label{Q-def}
Q^{k,m}_{n}(r) \ \equiv \ \frac{r^{m} P_{n}^{(k,m)}(z)}{\sqrt{N^{k,m}_{n}}} ; \  
N^{k,m}_{n} \equiv  \frac{W^{k,m}_n}{2^{2+k+m}}.
\Eeq
which implies an orthonormal basis 
\Beq
\int_{0}^{1} Q^{k,m}_{n}(r)Q^{k,m}_{n'}(r)\,(1-r^2)^k \,r \dd{r} \ = \ \delta_{n,n'}.
\Eeq
From now on, the $m$ index corresponds to the Fourier mode number in \eq{Fourier series}.  
Factorials replace Gamma functions in \eq{Q-def} for integer values of $k$. Each element of the radial basis satisfies the analytic property in \eq{analytic}. We could consider a non-integer shift in the starting value of $m$ in \eq{Q-def}; \eg this could correspond to starting with Worland polynomials rather than Zernike \cite{livermore_jones_worland_2007}. The following analysis would work with such an alteration, albeit with a large sacrifice to simplicity.  We therefore only discuss integer-$m$ values. 
 
The derivative operators for Jacobi polynomials imply the derivative operators on the radial basis
\Beq
\label{DR(r)+}
\left[\frac{d}{dr} - \frac{m}{r}\right]  Q^{k,m}_{n}(r) &=&  2 \sqrt{n(n+k+m+1)}\,Q_{n-1}^{k+1,m+1}(r),\n
\label{DR(r)-}
\left[\frac{d}{dr} + \frac{m}{r}\right]  Q^{k,m}_{n}(r) &=& 2\sqrt{(n+m)(n+k+1)} \,Q_{n}^{k+1,m-1}(r).
\Eeq
\Eqs{DR(r)+}{DR(r)-} form the foundation of a compact calculus of fields expanded in terms of $Q^{k,m}_{n}(r)$. 

Translating explicit statements like \eqs{DR(r)+}{DR(r)-} into more abstract operators acting on vector spaces allows easier interpretation of many of the results and derivations. 
Therefore, we use Dirac's bra-ket notation\footnote{see Appendix B for an introduction to the relevant aspects of this notation.} for the following analysis \cite{dirac_1939}. The following row-vector definition forms the  foundation for our algebra, 
\Beq
\bra{k,m,r}  \ \equiv  \ \left[ \, Q^{k,m}_{0}(r) , \, Q^{k,m}_{1}(r), \, \ldots \ \right].
\Eeq
The $k,m$ indices parameterise different bases. 
The corresponding ket is the transpose $\ket{k,m,r} = \bra{k,m,r}^{T}$. In the bra-ket notation \eqs{DR(r)+}{DR(r)-} become
\Beq
\label{basis-action}
\frac{1}{\sqrt{2}}\left[\frac{d}{dr} \mp  \frac{m}{r}\right] \bra{k,m,r}  \ \equiv \ \bra{k+1,m\pm1,r} D_{k,m}^{\pm}, 
\Eeq
where $D_{k,m}^{\pm}$ represent single-band matrices. The operators $D^{\pm}$ depend on the $k,m$ indices, but from now on we omit these labels when the context makes their values clear. $D^{-}$ is a diagonal matrix, and $D^{+}$ contains only a single band on the first super diagonal. The special case $m=0$ permits only the raising operation, and $D^{-} \equiv D^{+}$. In \eq{basis-action}, we choose the convention that $r$-space actions operate on the left, and matrix actions operate from the right. 

We see that the derivative operators on \eqs{m-up}{m-down} provide $m$-raising and lowering operators on the radial basis. \Eqs{DR(r)+}{DR(r)-} in many ways resemble the same operators applied to Bessel functions,
\Beq
\left[\frac{d}{dr} - \frac{m}{r}\right] \! J_{m}(r) =  -J_{m+1}(r),\quad \left[\frac{d}{dr} + \frac{m}{r}\right] \! J_{m}(r) =J_{m-1}(r).
\Eeq
This is no accident. Both Bessel functions and the scaled Jacobi polynomials behave similarly near the origin.  Alternatively, for many reasons Bessel functions are not as useful for representing the solutions to arbitrary PDEs \cite{boyd2011comparing}. In general, Bessel-function solutions achieve only algebraic convergence rates. For example, $f = r^{2} -1$ satisfies $\nabla^{2} f = 4, \ f(r=1) = 0$. This exact  polynomial solution only decays as $\mathcal{O}(n^{-5/2})$ in a Bessel-function series expansion. This behaviour results from a mismatch between the higher-order derivatives of an exact solution and a given approximating expansion; \eg the third derivative of $r^{2}-1$ vanishes identically, which does not happen for any finite Bessel function series. The same principle favours orthogonal polynomials over special functions in many other cases \cite{orszag_1971}.

\section{Tensor calculus}

We start with a scalar function 
\Beq
f(r,\theta) =  \sum_{m=-\infty}^\infty f_{m}(r)e^{im\theta}.
\Eeq
Using the bra-ket notation, we expand the radial dependence of the Fourier components in terms of the weighted Jacobi polynomials 
\Beq
f_{m}(r) \ = \ \braket{0,m,r}{f_{m}} \ = \ \sum_{n=0}^{\infty} f_{m,n}Q_{n}^{0,m}(r).
\Eeq
We expand all scalars, vector, and tensor components in terms of the $k=0$, and appropriate $m$ basis. As we differentiate a given component, the basis changes its $m$ index, and the $k$ index increments upward. A later section defines operators to convert between different $k,m$ bases. 

The covariant derivative does not couple the azimuthal Fourier modes. We can therefore write the action of its radial and azimuthal components on an individual Fourier mode as
\Beq
 \nabla_{\! r}f_{m}(r) \ = \ \frac{d}{dr}f_{m}(r),\quad \nabla_{\! \theta}f_{m}(r) \ = \ \frac{i m}{r}f_{m}(r).
\Eeq
Using \eqs{DR(r)+}{DR(r)-},
\Beq
 \nabla_{\!r} f_{m}(r) \ = \ \frac{1}{\sqrt{2}}\bra{1,m-1,r}D^{-} \ket{f_{m}} + \frac{1}{\sqrt{2}}\bra{1,m+1,r}D^{+} \ket{f_{m}}, \n
\nabla_{\!\theta} f_{m}(r) \ = \ \frac{i}{\sqrt{2}}\bra{1,m-1,r}D^{-} \ket{f_{m}} - \frac{i}{\sqrt{2}}\bra{1,m+1,r}D^{+} \ket{f_{m}}.
\Eeq

Introducing a spinor basis simplifies the computations of vectors and tensors. The spinor basis elements are
\Beq
e_{\pm} = \frac{1}{\sqrt{2}} \left( e_{r} \mp i e_{\theta} \right),
\Eeq
where $e_{r},e_{\theta}$ represent the unit vectors in the $r,\theta$ coordinate directions respectively.
A $2\times2$ unitary matrix transforms between the components of a vector in the $r,\theta$ basis and the $\pm$ basis.
In the spinor basis, the covariant derivative becomes
\Beq
\nabla \ = \ e_{r}\nabla_{\! r} + e_{\theta}\nabla_{\! \theta} \ = \ e_{-} \nabla_{\!-} + e_{+} \nabla_{\!+},
\Eeq
where acting on individual scalar functions, 
\Beq
\label{spin-grad}
\nabla_{\!\sigma} \ \equiv \ \frac{1}{\sqrt{2}}\left[\frac{d}{dr} - \frac{\sigma m}{r} \right],
\Eeq
with $\sigma = \pm1$. From standard vector calculus,  
\Beq
\nabla_{\!r}\,e_{r} =  \nabla_{\!r}\,e_{\theta} = 0,\quad \nabla_{\!\theta}\,e_{r} = \frac{1}{r}e_{\theta}, \quad \nabla_{\!\theta}\,e_{\theta} = -\frac{1}{r}e_{r}.
\Eeq
The spin basis diagonalises the connection coefficients such that
\Beq
\label{spin-connexion}
\nabla_{\!\sigma}\,e_{\mu} \ = \ - \frac{\sigma \mu}{\sqrt{2}r}e_{\mu},
\Eeq
where $\sigma, \mu = \pm1$. From the standard three-dimensional vector algebra of $e_{r},e_{\theta}, e_{3}$, we can easily deduce 
\Beq
e_{+}\times e_{-} \ = \ i\,  e_{3}, \quad e_{+}\cdot e_{-} \ = \ 1, \quad e_{+}\cdot e_{+} \ = \ e_{-}\cdot e_{-} \ = \ 0,
\Eeq 
where the dot product represents simple component-wise contraction, \ie \textsl{not} the proper complex-valued inner product.  Further, 
\Beq
\label{cross-prod}
e_{3}\times e_{\pm} \ = \ \pm  i\,  e_{\pm}.
\Eeq
We note that $e_{3}$ denotes any locally orthogonal third unit vector. For example, it could represent the axial direction in cylindrical coordinates, or an additional angle in a torus. 

\Eqs{spin-grad}{spin-connexion} imply the following formula for the covariant derivative of an arbitrary tensor component 
\Beq
\nabla_{\sigma} e_{\mu}\bra{k,m+\mu,r}  &=& - \frac{\sigma \mu}{\sqrt{2}r}e_{\mu}\bra{k,m+\mu,r} + \nonumber   e_{\mu}\nabla_{\sigma} \bra{k,m+\mu,r} \nonumber \\ &=& \frac{e_{\mu}}{\sqrt{2}}\left[\frac{d}{dr} - \frac{\sigma (m+\mu)}{r} \right]\bra{k,m+\mu,r} \nonumber \\ 
&=& e_{\mu}\bra{k+1,m+\mu+\sigma,r} D^{\sigma}.
\Eeq
Moreover, this formula holds if $\mu$ represents a multi-index such that
\Beq
e_{\mu} \ = \ \prod_{i=1}^{s} e_{\mu_{i}}, \quad \mathrm{and} \quad \bar{\mu} = \sum_{i=1}^{s} \mu_{i},
\Eeq
where $\mu_{i} = \pm1$ for each $i$.
If $\mu \in S = \left[-1,+1\right]^{s}$ represents an $s$-component multi-index, then 
\Beq
\mathrm{T}_{m}(r) \ = \ \sum_{\mu \in  S } e_{\mu}\braket{k,m+\bar{\mu},r}{ T_{m}^{\mu}} 
\Eeq
represents a given rank-$s$ tensor. 
The covariant derivative becomes the rank-$(s+1)$ tensor
\Beq
\label{tensor}
\nabla\, \mathrm{T}_{m}(r) \ = \ \sum_{\sigma = \pm1 } \sum_{\mu \in  S } e_{\sigma} e_{\mu}\bra{k+1,m+\bar{\mu}+\sigma,r}D^{\sigma} \ket{T_{m}^{\mu}} .
\Eeq
In \eq{tensor}, $k$ represents whichever basis the tensor began with. It need not equal the tensor rank, $s$. In most cases $k=0$.

Computing the gradient of a scalar 
\Beq
\label{grad}
\nabla f_{m}(r) \ = \ e_{-}\bra{1,m-1,r}D^{-} \ket{f_{m}} + e_{+} \bra{1,m+1,r}D^{+} \ket{f_{m}}.
\Eeq
We can further compute the second covariant derivative of a scalar
\Beq
\nabla \nabla f_{m}(r) &=&  e_{-} e_{-} \bra{2,m-2,r}D^{-} D^{-} \ket{f_{m}} +\nonumber \\ && e_{+} e_{-} \bra{2,m,r}D^{+} D^{-} \ket{f_{m}}  + \nonumber \\ && e_{-} e_{+} \bra{2,m,r}D^{-} D^{+} \ket{f_{m}} + \nonumber \\ && e_{+} e_{+} \bra{2,m+2,r}D^{+} D^{+}\ket{f_{m}}.
\Eeq

The operators $D^{\pm}$ commute in the sense that
\Beq
\label{D-commute}
D^{-}_{k+1,m+1} D^{+}_{k,m} \ = \ D^{+}_{k+1,m-1} D^{-}_{k,m},
\Eeq
where the $k,m$ indices show that operators take on different meaning depending on their order. 
This implies the covariant derivative components commute; as they must for a flat space. Taking the trace of the second derivative gives the Laplacian of a scalar
\Beq
\label{Laplacian}
\nabla^{2} f_{m}(r) &=&  2\bra{2,m,r}D^{+} D^{-} \ket{f_{m}}.
\Eeq

\section{Multiplication and Non-Constant Coefficients}

Many applications require multiplying scalar and vectors function by a non-constant, axisymmetric function of $r$. Here, we present a pair of operators representing multiplication without spoiling the bandedness of the matrix. In their discussion of vector formulations on disks, Sakai \& Redekopp (2009) used vector components of the form $r v_{m}(r)$ \cite{sakai2009application}. These components behave like scalars. Two additional Jacobi polynomial recursions reconcile this formulation with our current analysis. \Eq{a,b-1}, and \eqs{a+1,b}{exchange} together imply two distinct $r$-multiplication operators 
\Beq
\label{r+-}
r\, \bra{k,m,r} \ = \ \bra{k,m\pm1,r} R_{k,m}^{\pm},
\Eeq
where the dependency on $k$ and $m$ is dropped when inferred from context.
$R^{+}$ has entries on the diagonal and first super-diagonal, while $R^{-}$ contains entries on the diagonal and first sub-diagonal. For $m=0$, $R^{-} \equiv R^{+}$. \Eqs{R+ elements}{R- elements} show the individual elements associated with the operators in \eq{r+-}.

Given the gradient vector in \eq{grad},
\Beq
\label{rgrad}
r \nabla f_{m}(r) \ = \ e_{-}\bra{1,m,r}R^{+} D^{-}\ket{f_{m}} + e_{+} \bra{1,m,r}R^{-} D^{+} \ket{f_{m}}.
\Eeq
Both components of \eq{rgrad} contain the same $m$ dependence as the scalar $f_{m}(r)$. 

Considering \eq{z-def} and both $\pm$ parts of \eq{r+-},
\Beq
z \, \bra{k,m,r} \ = \ \bra{k,m,r} (2R^{-}R^{+} - I),
\Eeq
which leaves $k,m$ unchanged.
In the same style as \eq{D-commute}, radial multiplication commutes with itself,
\Beq
\label{R-commute}
R^{-}_{k,m+1} R^{+}_{k,m} \ = \ R^{+}_{k,m-1} R^{-}_{k,m},
\Eeq
where we note that the $\pm$ operators act on different $m$-bases depending on their order.

The standard three-term recursion relation for Jacobi polynomials implies the symmetric three-term recursion for the radial functions
\Beq
\label{recursion}
B^{k,m}_{n}\, Q_{n-1}^{k,m}(r) + A_{n}^{k,m}\, Q^{k,m}_{n}(r) + B_{n+1}^{k,m}\,Q^{k,m}_{n+1}(r) = z \, Q^{k,m}_{n}(r),
\Eeq
where
\Beq
\label{A-recur}
A_{n}^{k,m} &\equiv& \frac{m^2 - k^2}{(2n+k+m)(2n+k+m+2)}\n
\label{B-recur}
B^{k,m}_{n} &\equiv& \frac{2}{2n+k+m}\sqrt{\frac{n(n+k)(n+m)(n+k+m)}{(2n+k+m)^2 - 1}}.
\Eeq

Therefore, 
\Beq
Z  \ \equiv \ 2R^{-}R^{+} - I 
\Eeq
represents the standard symmetric tridiagonal Jacobi matrix for Jacobi polynomials. 
The operation of multiplication by any purely radial function occurs via 
\Beq
F(2r^2-1)\bra{k,m,r} \ = \ \bra{k,m,r} F(Z).
\Eeq 
We expand $F(z)$ in terms of a truncated series of any polynomials,
\Beq
\label{F-sum}
F(Z) \ = \ \sum_{n=0}^{N_F} F_{n} P_{n}(Z)
\Eeq
 An orthogonal polynomial basis allows the advantage of a matrix-valued three-term recursion of the form,
\Beq
\label{matrix-recursion}
B_{n+1} P_{n+1}(Z) + A_{n} P_{n}(Z) + B_{n}\, P_{n-1}(Z) = Z.P_{n}(Z) 
\Eeq
\Eq{matrix-recursion} need not represent any of the recursions in \eq{recursion}, but it can. For example, if we choose to expand in terms the same basis as the axisymmetric $k=m=0$ modes, then  
\Beq
A_{n} = 0, \quad B_{n} &\equiv&  \frac{n}{\sqrt{4n^2 - 1}}
\Eeq
But this is not an essential restriction. We could use Chebyshev polynomials for $P_{n}(Z)$ if the series in \eq{F-sum} converges more rapidly in that basis. In that case,
$A_{n} = 0$,  $B_{n} = 1/2$ for $n>0$, and $B_{0} = 1/\sqrt{2}$.

Reducing the truncation degree, $N_{F}$, reduces the bandwidth of the multiplication operator matrices and can save significant computation time even for modest differences in the truncation degree. $Z$ is a tridiagonal matrix, therefore the bandwidth of $F(Z)$ is $2N_{F}+1$, which is typically much smaller than the required size of the matrix.  

Given any choice in \eq{matrix-recursion} and taking $P_{-1} \equiv 0$, we can use Clenshaw's algorithm to compute $F(Z)$ directly via
\Beq
K_{N_{F}+1} = K_{N_{F}+2} = 0  
\Eeq
\Beq
K_{n} \ = \ F_{n} I + \frac{1}{B_{n+1}} (Z - A_{n} I). K_{n+1}  - \frac{B_{n+1}}{B_{n+2}}K_{n+2}
\Eeq
\Beq
F(Z) \ = \ P_{0}(Z) K_{0},
\Eeq
where $I$ represents the identity matrix, and $P_{0}(Z) = P_{0} I$ is a constant times the identity matrix. Each iteration of Clenshaw's algorithm increases the bandwidth by one, via multiplication by the tridiagonal matrix $Z$. If the coefficients vary strongly with latitude, then the multplication matrices are dense. But this happens in a gradual way. Most of the time, the coefficients do not vary nearly as rapidly as the solution itself. 

Multiplication by non-azimuthally symmetric functions must obey more complicated selection rules for the different $m$  indices.  This would create systems that couple both the $\theta$ and $r$ directions. 
We do not consider linear operators with $m>0$ dependence.

\section{Basis conversion}

\Eq{Laplacian} has a problem: the Laplacian of a scalar with $k=0$ results in a scalar with $k=2$. We must convert between the $k$ bases in order to compare the input and output.

\Eq{a+1,b} implies  a 2-band conversion operator such that 
\Beq
\label{conversion C}
 \bra{k,m,r}  \ \equiv \ \bra{k+1,m,r} C_{k,m},
\Eeq
where the dependence on $k$ and $m$ is again dropped when inferred from context. $C$ has entries on the diagonal and first super-diagonal. Like the $D^{\pm}$ operators, we omit the $k,m$ indices when the context makes it clear. \Eq{C elements} shows the individual elements associated with the operators in \eq{conversion C}.

Raising the $k$ index allows converting to the basis that most naturally represents derivatives. This is the same principle that underlies sparse representations with Chebyshev polynomials, \ie
\Beq
\frac{d}{dz}T_{n}(z) = n U_{n-1}(z), \quad 2 T_{n}(z)  = U_{n}(z) - U_{n-2}(z). 
\Eeq
Many papers reinventing Chebyshev numerical schemes rely on this simple fact \cite{zebib_1984,coutsias_hagstrom_1996,boyd_book,greengard_1991,doha2002efficient,julien_watson_2009,doha2009efficient,muite_2010,olver2013fast,viswanath_2015,townsend_olver_2015}.  Chebyshev polynomials are special cases of Jacobi polynomials with $a=b=-1/2$ and $a=b=1/2$ representing the first and second kind, respectively. Using greater flexibility in the $a,b$ indices allows for an increased ability to construct sparse operators for a particular application.

\subsection{Heisenberg Algebra}

Radial multiplication, differentiation, and basis conversion satisfy simple algebraic compatibility conditions. Differentiation and multiplication satisfy standard commutation relations, but with the added element of conversion. 

\begin{figure}
\begin{center}
\includegraphics[scale=0.35]{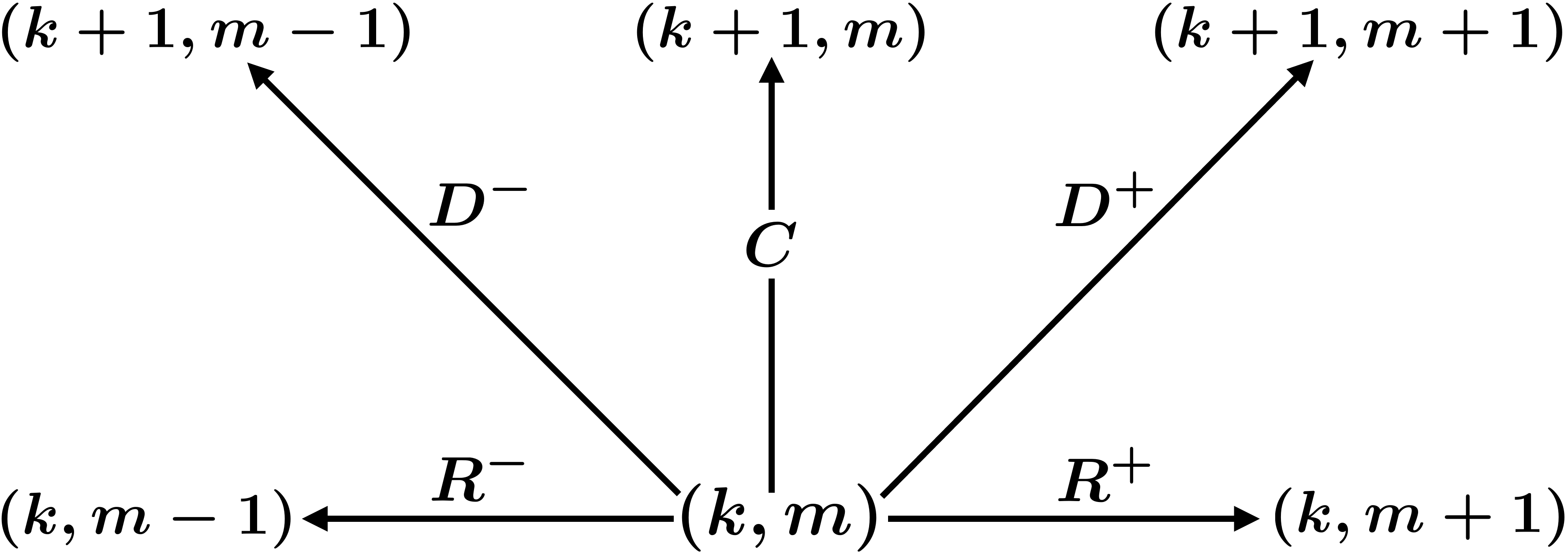}
\caption{Diagram of operator structure for the $(k,m)$ indices. Every other node in the graph contains its own set of operators mapping to new nodes in a similar fashion.  Only the right half of the diagram exists for $m=0$. 
\label{ops-structure}}
\end{center}
\end{figure}

The commutation relations satisfy the two-dimensional Heisenberg Lie algebra, $\mathfrak{h}_{2}$. That is,
\Beq
\label{H-algebra}
\left[ D^{+},R^{-} \right] \ = \ \left[ D^{-},R^{+} \right]  \ = \   \sqrt{2} C.
\Eeq
Keeping track of the different $k,m$ indices,
\Beq
D_{k,m\mp1}^{\pm}R_{k,m}^{\mp} - R_{k+1,m\pm1}^{\mp}D_{k,m}^{\pm}  \ = \ \sqrt{2} C_{k,m},
\Eeq 
Otherwise all other $\pm$ operators commute with each other and with $C$. 
\Fig{ops-structure}{} shows the structure of the different elements of $\mathfrak{h}_{2}$ acting on the $\bra{k,m,r}$ state. The diagram omits the equivalent operators acting on the other nodes in the infinite $(k,m)$ lattice. Traditional representation theory considers the Lie algebra of operators acting on a single Hilbert space. Here it appears that we are considering the unconventional situation of a Lie algebra of operators acting \textsl{between} Hilbert spaces; each with different $k,m$. Considering a total Hilbert space comprising the direct sum of all individual $k,m$ spaces rectifies this issue. 

The operators also satisfy conditions apart from commutation relations. Taking the $\theta$ component of \eq{rgrad} implies
\Beq
m C  =\frac{1}{\sqrt{2}}\left( R^{+}D^{-} - R^{-} D^{+} \right).
\Eeq
Taking the $r$ component of \eq{rgrad} implies,
\Beq
\label{rd/dr}
 r \frac{d}{dr} C = \frac{1}{\sqrt{2}}\left ( R^{+}D^{-} + R^{-}D^{+} \right).
\Eeq
\Eq{rd/dr} is equivalent to equation (15) from Matsushima \& Marcus (1995) \cite{matsushima1995spectral}. 

The Heisenberg algebra allows the automatic reformulation of equations at the numerical level. \Eq{H-algebra} gives a simple and reliable recipe for commuting radial non-constant coefficients past derivatives using the numerical operators. 

\subsection{Boundary restriction}

Solving PDEs on the unit disk usually requires appropriate boundary conditions at $r=1$. No boundary at $r=0$ is needed, nor allowed.  To impose the boundary condition at $r=1$ we use  
\Beq
\label{boundary-restriction}
Q_{n}^{k,m}(r=1) \ = \ \sqrt{2(2n+m+k+1){{n+k}\choose{k}}{{n+m+k}\choose{k}}},
\Eeq
which follows from \eqs{P at z=1}{Q-def}.
We can build higher-order boundary restriction operators out of \eq{boundary-restriction} and various  differentiation formulae. For example,
\Beq
\label{d/dr at r=1}
\frac{d}{dr}Q_{n}^{k,m}\Big|_{r=1} =  \frac{2n(n+k+1)+ m(2n+k+1) }{k+1} Q_{n}^{k,m}(r=1).
\Eeq

It is convenient to modify the basis so that the boundary condition only appears in the first element.   This leads to banded systems, as opposed to almost banded systems, which may have dense rows corresponding to the boundary conditions.   Many black-box linear algebra packages can solve banded systems very efficiently, while almost-banded systems require specially designed schemes, \eg  the adaptive QR method \cite{olver2013fast}, or use of the Woodbury formula.

We describe how to modify the the basis to accommodate Dirichlet boundary conditions, corresponding to \eqref{z=1}.  The adjoint of the conversion operators gives a $k$-lowering operator.
\Beq
\label{C dag operator}
(1-r^2) \bra{k,m,r} \ = \ \bra{k-1,m,r} C^{\dag}
\Eeq
where $C^{\dag}$ has entries on the diagonal and first sub-diagonal. This follows from a straightforward application of \eqs{a,b-1}{exchange}. \Eq{Cdag elements} gives the explicit indices of \eq{C dag operator}.
The $k$-lowering operator satisfies the commutation relations,
\Beq
[ C^{\dag},D^{\pm}] \ = \ \sqrt{2} R^{\pm}.
\Eeq
Otherwise $C^{\dag}$ commutes with all other relevant operators. The $R^{\pm}$ operators are the adjoints of each other. Adjoints of the $D^{\pm}$ operators exist, but provide little apparent use for the types of problems we consider here.  

Define
\Beq
\label{B-op}
B \ \equiv \ C^{\dag} +  \alpha^{-1} \ket{0}\bra{0},\quad
\mathrm{where}\quad
\alpha \equiv \braket{k,m,r=1}{0},
\Eeq
where $\bra{0}$ represents the unit basis element corresponding to the $n=0$ mode.
\Eq{B-op} defines a two-band lower-diagonal operator. 

As before,
\Beq
f_{m}(r) \ = \ \braket{k,m,r}{f_{m}}
\Eeq
If we change variables such that,
\Beq
\ket{f_{m}} \ \equiv \ B \ket{g_{m}}
\Eeq
then
\Beq
f_{m}(r) & = &  \bra{k,m,r}C^{\dag} \ket{g_{m}} + \frac{\braket{k,m,r}{0} }{\braket{k,m,1}{0}} \, \braket{0}{g_{m}} \nonumber \n &= & (1-r^2)\braket{k+1,m,r}{g_{m}} + \frac{\braket{k,m,r}{0} }{\braket{k,m,1}{0}} \, \braket{0}{g_{m}}.
\Eeq
Therefore
\Beq
f_{m}(r=1) \ = \ \braket{0}{g_{m}}.
\Eeq
The $B$ operator localises the boundary-restriction operator to the first mode of $g$. This contrasts with a dense row operation for the original variable. This comes at the expense of a single band in the bulk equation operators, but produces an entirely banded system.

For the derivative, because for every individual element $\bra{k,m,r=1} > 0$, there exists a diagonal operator $\Lambda$ such that  
\Beq
f_{m}'(r=1) &=& \bra{k,m,1}\Lambda \ket{f_{m}} 
\Eeq
This also extends to more general left-hand side restriction operators. In the case of the derivative, \eq{d/dr at r=1} implies the diagonal elements
\Beq
\label{Lambda_n}
\bra{n} \Lambda \ket{n} \ = \ \frac{2n(n+k+1)+ m(2n+k+1) }{k+1}.
\Eeq
Therefore, if
\Beq
\label{general basis-recombination}
\ket{f_{m}} \ \equiv \ \Lambda^{-1} B \ket{g_{m}},
\Eeq
then 
\Beq
f_{m}'(r=1) \ = \ \bra{k,m,1} B \ket{g_{m}} \ = \ \braket{0}{g_{m}}.
\Eeq
The same principle applies for any $\Lambda$, not only for the derivative. Julien \& Watson (2009) demonstrate many of the advantages of a sparse transformation to a Galerkin basis for general boundary operators \cite{julien_watson_2009}.
\Eq{Lambda_n} shows that $\Lambda_{n=0} = 0$ for $m=0$, otherwise $\Lambda_{n} > 0$. This precludes forming \eq{general basis-recombination} for $n=m=0$. Neumann boundary conditions fail to determine this particular mode, which requires additional information.  In many applications, setting this mode requires fixing a physically irrelevant (but nevertheless mathematically essential) gauge condition.

\subsection{Quadrature}

All the above analysis considers field variables in the spectral coefficient domain. This allows for efficient differentiation operations. Given a set of spectral coefficients, we need the option to construct the function on a radial grid. We also need to compute spectral coefficients from data on a grid.  In general, we need each $m$ mode on the same radial grid. This allows constructing solutions on the full $r,\theta$ disk.  We consider an expansion in terms of a finite number of radial modes,
\Beq
f_{m}(r) \ = \ \sum_{n=0}^{N_{k,m}} Q_{n}^{k,m}(r) \braket{k,m,n}{f_{m}},
\Eeq
where $N_{k,m}$ depends on both $k$ \& $m$. We use Gauss-Legendre quadrature on the grid $z_{i} = 2r_{i}^{2}-1$ such that
\Beq
\label{quadrature}
\braket{k,m,n}{f_{m}} = \sum_{i=1}^{N_{r}} Q_{n}^{k,m}(r_{i}) (1-r_{i}^2)^{k} w_{i} f_{m}(r_{i})
\Eeq
for some $N_{r}$. Gauss quadrature is exact for all $z$ polynomials less than degree $2N_{r}$, or $r$ polynomials less that degree $4N_{r}$. In the $r$ variable,
\Beq
\deg \left[Q_{n}^{k,m}(r)(1-r^2)^{k}f_{m}(r) \right] = 2 (n +m+ k+N_{k,m}) < 4 N_{r}
\Eeq
where $n \le N_{k,m}$. Therefore,
\Beq
\label{Nkm}
N_{k,m} \ \equiv \ N_{r} -1 - \mathrm{floor}\!\left[{\frac{k+m}{2}}\right] 
\Eeq
yields exact quadrature on a Gauss-Legendre grid. Gauss-Legendre quadrature is based on the $k=0,m=0$ case. For non-integer starting $k$ values, $w_{i},r_{i}$ would result from a more general Gauss-Jacobi quadrature, but the basic principle remains the same. In practice, we never need to compute \eq{quadrature} with $k > 0$. Instead, we can project any function on $r_{i}$ to the $k=0$ spectral coefficients and convert the result up the desired $k$ using the $C$ conversion operator. 

The absence of a fast transform is the largest downside to the Jacobi-based method compared to a Chebyshev scheme. In extreme cases of large $N_r$, where a fast transform would show significant gains, Chebyshev schemes struggle to represent the severe analytic behaviour of the solution near $r=0$. For example, if $m$ is an even integer, with $p \equiv m/2$, and $k=0$, then
\Beq
\label{Chebyshev-series}
\nonumber && r^{m}P_{n}^{0,m}(2r^{2}-1) \ = \ \\
&&2\sideset{}{'}\sum_{j=0}^{n+p}T_{2j}(r)\sum_{i=\max\{j,p\}}^{n+p}
\frac{\binom{2 i}{i+j}
   \binom{i}{p} \binom{n+p}{i}
   \binom{i+n+p}{n}}{\binom{n+p}{
   p}}\frac{(-1)^{i + n + p}}{4^i}.
\Eeq
where $T_{2j}(r)$ represents a Chebyshev polynomial of even degree $2j$. The primed-summation notation in \eq{Chebyshev-series} represents multiplying the $j=0$ term of the series by $1/2$. 
The series in \eq{Chebyshev-series} contains Chebyshev modes up to degree $2n+m$. A similar relation exists for odd $m$ values. Given a relative error tolerance $\varepsilon$, in some cases  these series coefficient will become smaller than the allowed error before the end of the exact series. However, even this does not gain much advantage. In the regime $m \sim n$, a Chebyshev series expansion of $Q_{n}^{0,m}(r)$ requires nearly all of the $2n+m$ modes to achieve an accurate representation.  We note that these considerations ignore the possibility of a low-rank approximation to \textsl{specific} functions, which can also exploit a fast transform; see \cite{townsend_etal_2016}.

Also, because $N_{k,m}$ decreases with  larger $m$, a Jacobi-based scheme requires fewer total modes to obtain uniform resolution over the entire disk. That is, after summing \eq{Nkm} over $|m|\le N_{\theta}/2$ we find 
\Beq
N_{\mathrm{tot.}} \approx \left(N_{r}-\frac{N_{\theta}}{8}\right) N_{\theta},
\Eeq
compared to simply $\approx N_{r}N_{\theta}$ for a Chebyshev radial scheme. The savings are maximised for $N_{\theta} \approx 4 N_{r}$, with $N_{\mathrm{tot.}} \approx  N_{r}N_{\theta}/2$.

Two final mitigating factors arise from (\textit{i}) the somewhat-large cost of fast transforms for moderate sizes, and (\textit{ii}) the easy parallelisation of matrix transforms (\ie quadrature). Table~\ref{table:timings} compares the timing of a fast Chebyshev transform (FCT, based on the type-II fast cosine transform), versus quadrature computed with a matrix-vector multiply (MMT). The tests show that while the FCT gives better asymptotic scaling, the MMT does better at moderate sizes. The speed crossover occurs somewhere between $64$ and $128$ modes; a reasonable size for many applications.  The timing difference remains less than an order of magnitude all the way up to $1024$ modes. Lastly, a quadrature-based transform uses the orthogonality of the eigenfunctions to compute the inverse transform from the transpose of the forward transform. That is, both the forward and inverse transform represent $N_{r} N_{k,m}$ operations.  However, because of the matrix-vector structure, the operations actually entail $N_{r}$ totally independent vector-vector dot products. This parallelises quite easily. 

\begin{table}[ht]
\smallskip
\centering
\begin{tabular}{ccccccccccc}
\hline
$N_{r}$ & $:$ &               32 & 64 & 128 & 256 & 512 & 1024 & 2048  \\
\hline
FCT ($\mu$sec) & $:$ & 0.94 & 1.97 & 4.42 & 8.56 & 17.8 & 40.7 & 81.3 \\
\hline
MMT ($\mu$sec) & $:$ & 0.43 & 1.75 & 9.19 & 30.1 & 117. & 462. & 1905. \\
\hline
MMT/FCT & $:$ & 0.45 & 0.89 & 2.08 & 3.51 & 6.59 & 11.3 & 23.4 \\
\hline
\end{tabular}
\caption{Timings for Fast Chebyshev/Cosine Transform (FCT) versus a Gauss-Legendre Matrix-Multiply Transform for different resolutions. The timings use SciPy (FCT) and NumPy (MMT) on a Late-2013 Mac Pro 2.7 GHz 12-Core Intel Xeon E5. Each timing represents the average over 30,000 forward and inverse transforms of random data; only the actual transforms factor into the average.}
\label{table:timings}
\end{table}

\section{Examples \label{examples}}

This section considers a series of practical examples. Example 1 demonstrates the operator methods on a simple scalar equation. Example 2 demonstrates the same for a vector set of equations. Example 3 concerns non-constant coefficients in a vector system. Example 4 considers a problem with full $\theta,r$-dependence (\ie with a large number of $m$ values) and considers computational performance.  

The first three examples entail finding the generalised eigenvalues and eigenvectors for an analytically tractable, and/or previously published numerical test.  We compute eigenvalues because of the demanding nature of these problems. Eigenvalues and eigenvectors give a significant fraction of all the useful information about a given system. For example, the eigenvalue problem determines the stability of general time-stepping schemes. In some cases, defective numerical schemes can produce large spurious eigenvalues in the wrong part of the complex plane \cite{zebib_1984,gardner_etal_1989,mcfadden_etal_1990}. Such methods can give accurate solutions for a single linear solve with specific right-hand side terms. But rogue eigenvalues in linear operators can render a time stepper useless after a small number of iterations. A small random projection onto a bad eigenvector will diverge exponentially upon iteration. Also, eigenvalues in the wrong part of the complex plane can prevent otherwise conserved quantities from remaining constant. 

In each case we build the operator matrices by inspection of the original equation and use standard linear algebra libraries for the actual calculations. Apart from example 4, we code all the examples in Python, and take advantage of the features within NumPy and/or SciPy.  We produce Jupyter notebooks, and scripts for each case and make them available online within the Dedalus project\footnote{\url{dedalus-project.org}}.

Example 4 is coded in the ApproxFun.jl  Julia package \cite{olver2013fast,olver2014practical}, to utilize its implementation of the adaptive QR method.

\paragraph{\textbf{Example 1} (Bessel function eigenvalue problem):} We solve the following eigenvalue problem  
\Beq
\label{BesselJ}
\frac{1}{r}\frac{d}{dr} \left[r \frac{d}{dr} J_{m}(\kappa r)   \right] - \frac{m^{2}}{r^{2}} J_{m}(\kappa r) = - \kappa^{2}J_{m}(\kappa r),
\Eeq
with the Dirichlet boundary condition, 
\Beq
\label{BesselZero}
J_{m}(\kappa)=0.
\Eeq
Many available software packages compute Bessel functions and their corresponding zeros. We therefore use this example as a numerical test.  

We solve the generalised eigenvalue problem,
\Beq
2 D^{-}D^{+}\ket{f_{m}} = -\kappa^{2}\, CC \ket{f_{m}},
\Eeq
where $CC$ denotes applying two (different) conversion operators. The left-hand side operator only contains a single off-diagonal band. Each $C$ operator only contains a single off-diagonal band. The product of the two is tridiagonal. This contrasts with the penta-diagonal system for this problem in \cite{matsushima1995spectral}.

The boundary restriction operator in \eq{boundary-restriction} replaces the last matrix row on the left-hand side. A row of zeros replace the last matrix row on the right-hand side. After constructing the appropriate right- and left-hand side matrices, we use the NumPy eigenvalue solver to find the eigenvalues and eigenvectors in coefficient space.  We sort the results from lowest-to-highest in eigenvalue, and transform the coefficient-space vector to a radial grid; see \S5.3 for details regarding the spectral-to-spatial transforms. 

\begin{figure}
\begin{center}
\includegraphics[scale=0.35]{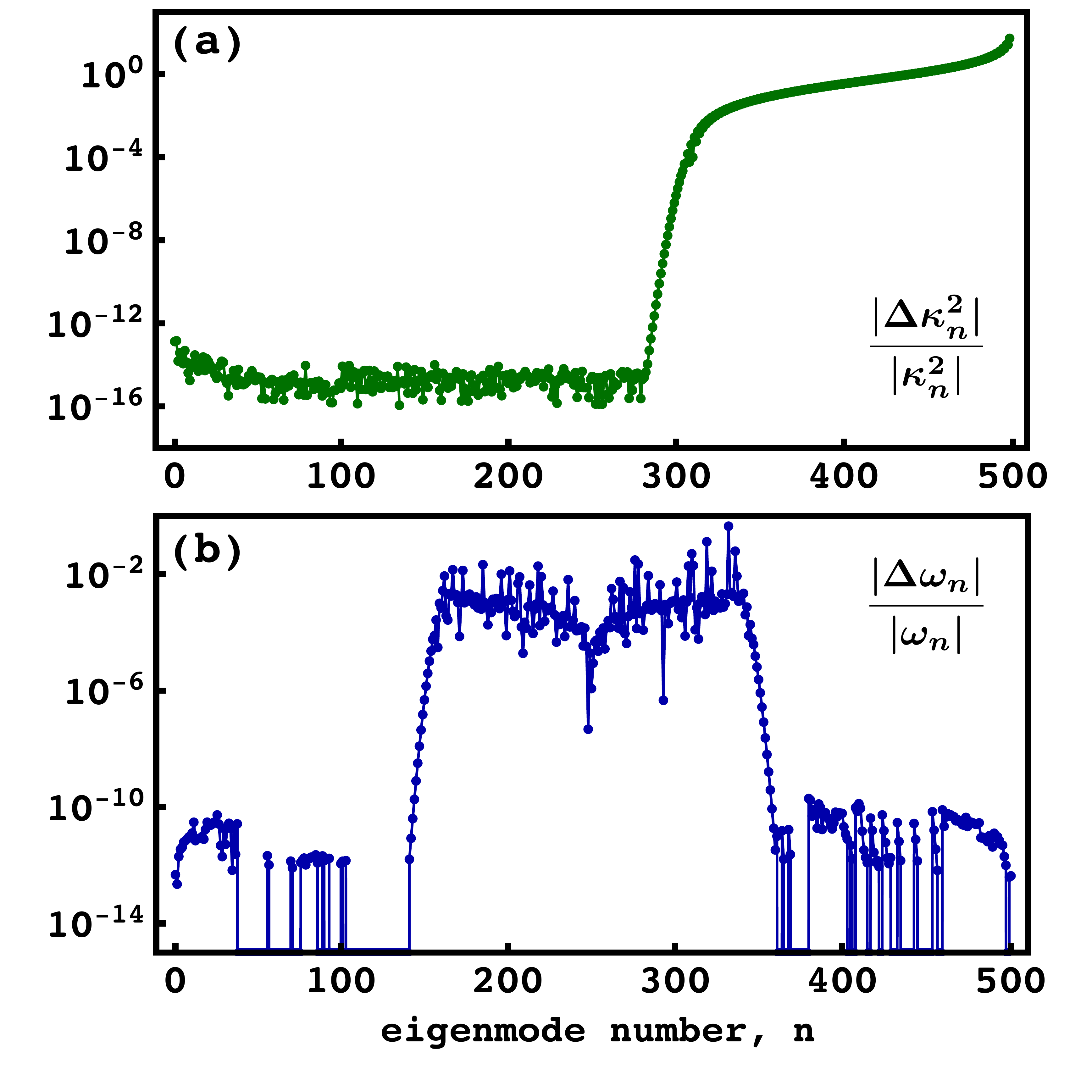}
\vspace{-0.5cm}
\caption{Plots of relative eigenvalue error from examples 1 \& 2. Each plot compares the numerical eigenvalue solution with the analytically computed value for 500 individual sorted modes.  Panel (a) corresponds to Bessel function example~1 for $ m=50$. Panel (b) corresponds to the inertial waves example 2 $ m=1$. Each case contains roughly $70\%$ "good" values with high accuracy, and $30 \%$ "bad" values with significant error. In both cases, the bad values still lie in the same region on the complex plane as the true solutions. The bad values correspond to small-scale modes that outstrip the fixed resolution. In the Bessel function example, small-scale modes correspond to large $\kappa$. Hence the error become large for large mode numbers. For inertial waves, small-scale modes correspond to slow frequencies with $| \omega | \approx 0$. The more accurate large-scale modes correspond to $|\omega| \approx 1$. Therefore, the significant errors on panel (b) occur for intermediate mode numbers; see Table (\ref{table:Eigenvalues})
\label{eigenvalue-plot}}
\end{center}
\end{figure}

Figure~\ref{eigenvalue-plot}a shows the relative error in the eigenvalues for an $m=50$ solution set with $N_r=500$ radial modes. The plot is typical of the behaviour for other $m$ values. We choose this case to examine closely because it gives a good illustration of the ability of the method to capture the large dynamic range near $r=0$. However the method's performance near the origin is independent of the $m$ value. The behaviour of the numerical eigenvalues is typical of orthogonal polynomial solutions of boundary-value problems, including those using Chebyshev polynomials. Typically, given a value of $N_{r}$, approximately the first $1/2$ to $2/3$ of the spectrum will contain values with virtually no error; apart from machine floating-point error, and a little from the eigenvalue solver routine. In figure~\ref{eigenvalue-plot}a this region spans from mode 0 to about mode $300$ out of $500$ modes. After this point, the exact solutions begin to outstrip the inherent finite-degree polynomial approximation. This error increases rapidly, with the largest eigenvalues scaling such that 
\Beq
|\Delta \kappa_{n}| \sim \kappa_{n}^{2}.
\Eeq
Like Chebyshev schemes, this behaviour is typical for orthogonal polynomial solutions of two-point boundary eigenvalue value problems. Figure~\ref{eigenvalue-plot} implies that obtaining accurate solutions to Poisson's equation,  
\Beq
\nabla^{2}f_{m}(r)  = s_{m}(r)\quad  \longleftrightarrow \quad 2 D^{-}D^{+} \ket{f_{m}} = CC \ket{s_{m}}
\Eeq
requires that $s_{m}(r)$ does not project onto the inaccurate part of the spectrum for a given resolution. This implies that the resolution needed to compute $f_{m}(r)$ is roughly $3/2$ to $2$ times larger than needed to resolve $s_{m}(r)$. But after reaching the critical number of modes, the error is effectively machine precision. This same principle holds for most types of equations. 

\begin{figure}
\begin{center}
\includegraphics[scale=0.36]{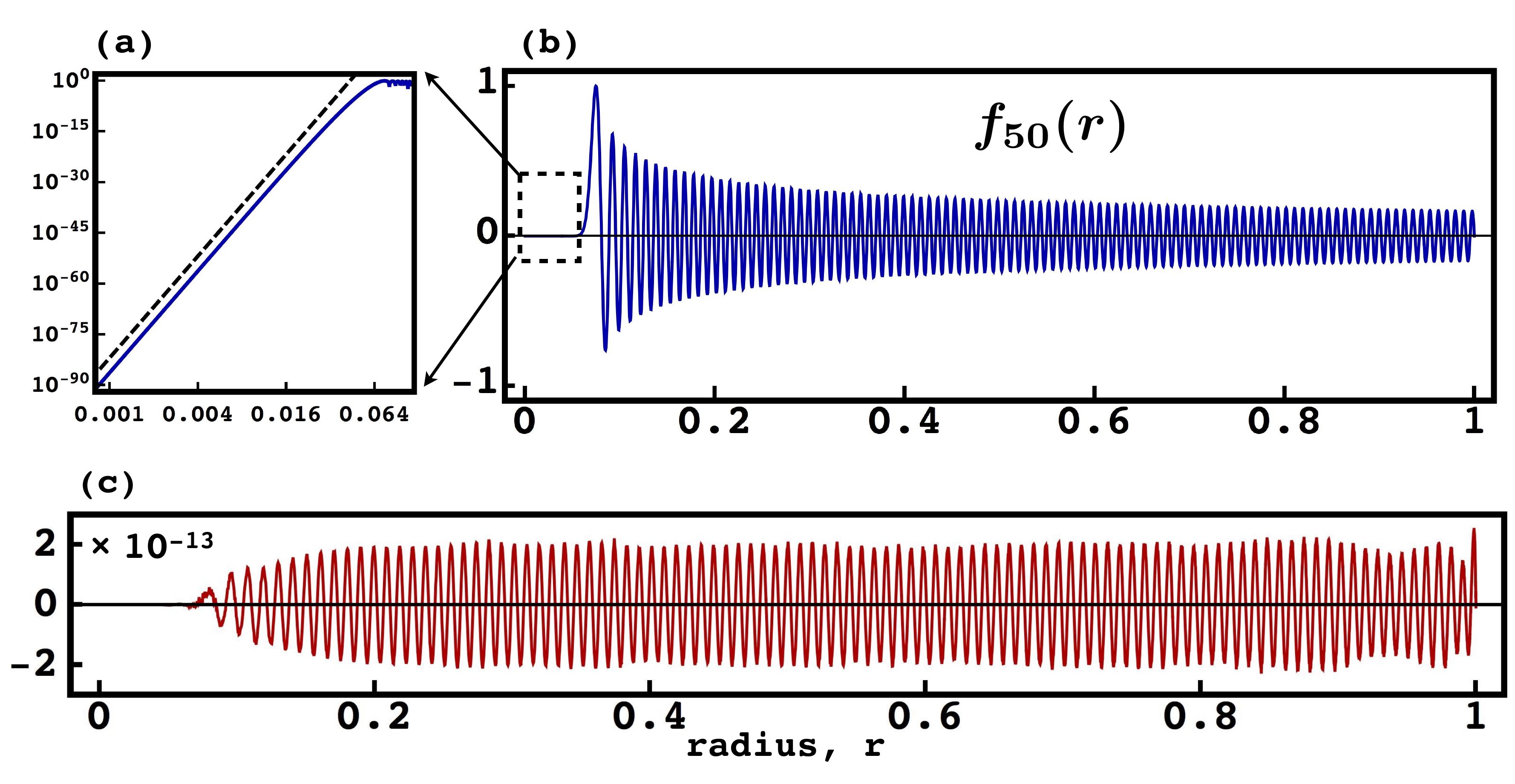}
\vspace{-0.3cm}
\caption{Plots for the $200$th eigenfunction from Example 1 for $m=50$, and $N_{r} = 500$. Panel (a) shows a $\log$-$\log$ plot of radial dependence near the origin. The solid blue line shows the computed eigenfunction. The eigenfunction should scale as $\sim r^{50}$. The dashed black line denotes this scaling.  Panel (b) shows the eigenfunction across the entire radius. The solution vanishes at $r=1$.  Panel (c) shows the error between the computed solution and an analytically computed Bessel function. The characteristic error is $\approx 2 \times 10^{-13}$, and mostly uniform over the domain. 
\label{Bessel-plot}}
\end{center}
\end{figure}

\begin{figure}
\begin{center}
\includegraphics[scale=0.40]{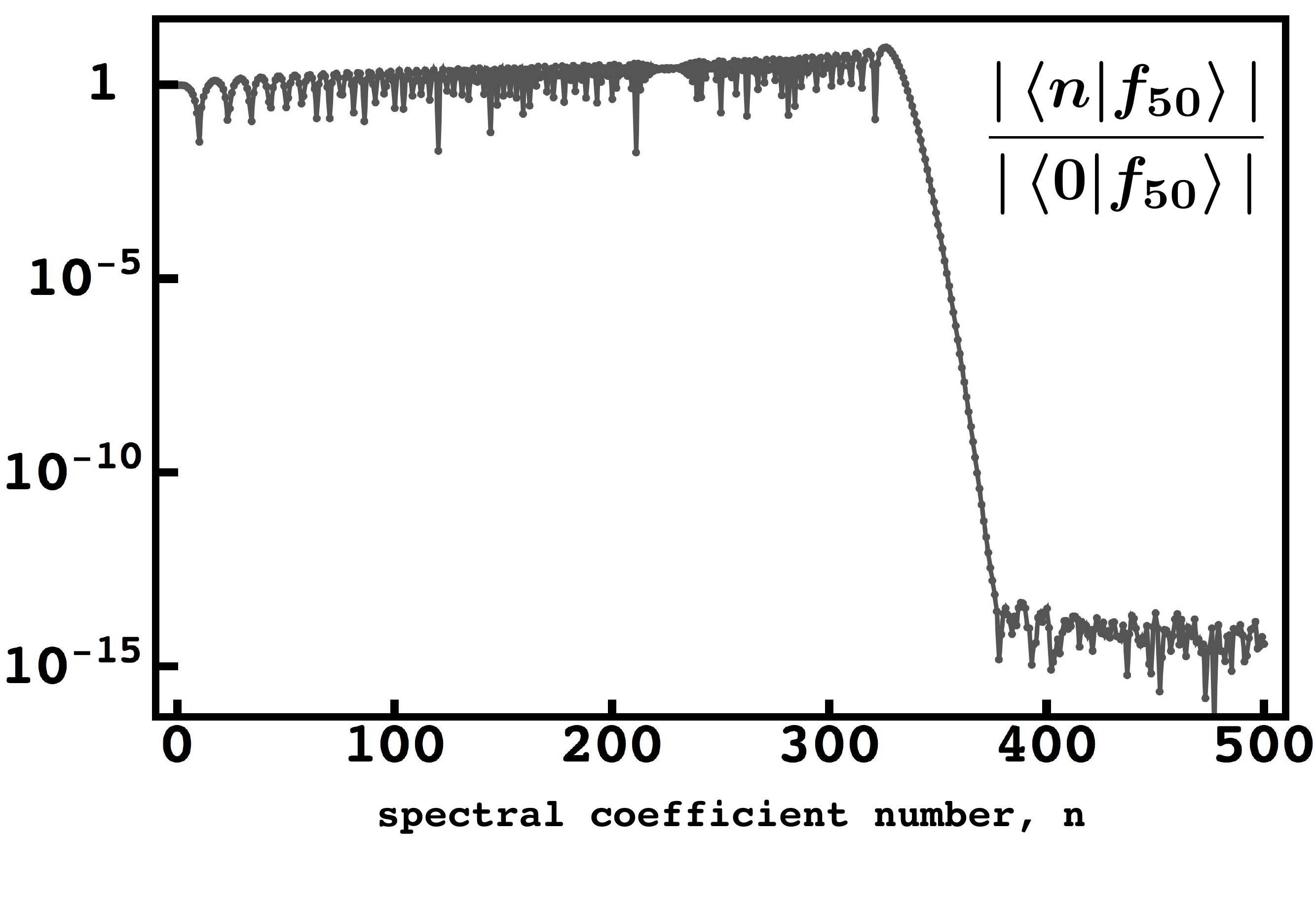}
\vspace{-0.5cm}
\caption{Plot of the spectral coefficient amplitude of the $200$th eigenfunction from Example~1 for $m=50$. In general, the $n$th eigenfunction requires $\approx\!2n$ spectral coefficients to resolve.  
\label{coef-plot}}
\end{center}
\end{figure}

Figure~\ref{Bessel-plot} shows the eigenfunction behaviour for the $n=200$ eigenvector for $m=50$. This corresponds to $\kappa_{200} \approx 707.447066905$. The highly oscillatory solution is accurate to roughly 13 digits mostly uniformly throughout the domain. Figure~\ref{coef-plot} shows the spectral coefficients for the same eigenfunction. The magnitude of the coefficients decrease dramatically after approximately $300$ modes. This accords with the similar behaviour in the the eigenvalues. Trying to represent this mode with fewer than $300$ to $400$ modes would incur significant error. Beyond this number, more resolution contributes little to the fidelity of the solution. 

\begin{table}[ht]
\smallskip
\centering
\begin{tabular}{cccccccc}
\hline
$n$ & $:$ & 0 & 100 & 200 & 300 & 400 & 500 \\
\hline
$\kappa_{n}\, (m=50)$ & $:$ & 5.71e1 & 3.91e2 & 7.07e2 & 1.02e3& 1.33e3& 1.64e3 \\
\hline
$\omega_{n}\, (m=1)$ & $:$ & -2.14e-1 & -3.17e-3 & -1.23e-3 & 1.20e-3& 3.12e-3& 3.19e-1 \\
\hline
\end{tabular}
\caption{Approximate analytically computed eigenvalues versus mode number. The approximate values are intended to give scale to the bottom axes on Figure~\ref{eigenvalue-plot}. $\kappa_{n}$ increases monotonically from large-scale modes to small scale. For $\omega_{n}$, small-scale modes correspond to the middle of the sorted eigenvalue spectrum, and large-scale modes correspond to the values closer to unity. Inertial waves break reflection symmetry. Hence the lack of symmetry $\omega \to -\omega$.}
\label{table:Eigenvalues}
\end{table}

Finally, we note that the matrices for radial component in polar coordinates are sparser than the Cartesian Chebyshev analogue in one dimension.  Chebyshev discretisation applied in polar coordinates must account for the geometric terms with non-constant coefficients. Using operators that conform to the disk geometry renders systems such as \eq{BesselJ} effectively constant coefficient problems.

\paragraph{\textbf{Example 2} (Inviscid inertial waves in a upright rotating cylinder):} The following example highlights vector-scalar coupled system with a non-trivial analytical solution. It therefore makes a good system for numerical comparison. We consider 
\Beq
\label{v-eq}
i \omega\,  v_{m}(r)  + e_{3} \times v_{m}(r) +  \nabla p_{m}(r) \ = \ 0, \n
\label{p-eq}
i \omega \nabla \cdot v_{m}(r) + \alpha^{2} p_{m}(r) \ = \ 0.
\Eeq
The two-dimensional velocity vector and pressure, $v_{m}(r) = v^{+}_{m}(r) e_{+} + v^{-}_{m}(r) e_{-}$ and  $p_{m}(r)$ represent a full eigenstate. The parameter  $\omega$ represents a frequency eigenvalue. The parameter $\alpha$ represents a mode aspect ratio in the third dimension. 
We solve \eqs{v-eq}{p-eq} with the normal-velocity boundary condition 
\Beq
\label{v_r=0}
e_{r} \cdot v_{m}(r=1) = 0.
\Eeq
Greenspan (1968) \cite{greenspan_1968} discusses at length the derivation and physical interpretation of \eqs{v-eq}{p-eq}. We take the system as given.  Inertial waves make a good test problem because of the non-trivial frequency constraint 
\Beq
-1 < \omega < 1.
\Eeq
As with the Bessel function problem, numerical eigenvalue solvers often produce a subset of values that differ significantly from their true solution. It is, however, a sign of a good scheme if all values are maintained in the same region of the complex plane as the true result. As discussed above, a poorly designed scheme can produce very large spurious eigenvalues with a different complex character. 

An analytical solution to \eqss{v-eq}{v_r=0} results from collapsing the vector system to the scalar problem 
\Beq
\label{analytical-inertial}
\nabla^{2} p_{m}(r)  \ = \ -\frac{(1-\omega^{2})\alpha^{2}}{\omega^{2}} p_{m}(r)
 \Eeq
with
\Beq 
\label{inertial-bc}
\omega \frac{d}{dr} p_{m}(r) + \frac{m}{r} p_{m}(r) \ = \ 0 \quad \mathrm{at} \quad r =1.
\Eeq
\Eq{analytical-inertial} requires a solution in terms of Bessel functions, $p_{m}(r) = J_{m}(\kappa r)$. \Eq{inertial-bc} requires
\Beq
\label{kappa-eq}
\kappa^{2} \ = \ \frac{(1-\omega^{2})\alpha^{2}}{\omega^{2}}, \quad \kappa \omega J_{m}'(\kappa) + m J_{m}(\kappa) \ = \ 0.
\Eeq
We solve \eq{kappa-eq} with a standard numerical root finder. 

Recall \eq{cross-prod}, $ e_{3} \times e_{\pm} = \pm i \, e_{\pm} $. This implies the following block eigensystem 
\Beq
\label{block-system}
\left[\begin{array}{ccc}  C & 0 & D^{+} \\ 0 & -C & D^{-} \\ 0 & 0 & \alpha^{2}C  \end{array}\right]  \ket{\Xi_{m}} = -\omega \left[\begin{array}{ccc} C & 0 & 0 \\ 0 & C & 0 \\ D^{-} & D^{+} & 0  \end{array}\right] \ket{\Xi_{m}},
\Eeq
where the eigenstate is the vectorisation of the velocity and pressure,
\Beq
\ket{\Xi_{m}} \equiv \left[\begin{array}{ccc} \left| iv_{m}^{+}\right> & \left|{iv_{m}^{-}}\right> & \left|{p_{m}}\right> \end{array}\right]^{T}.
\Eeq
Note that absorbing a factor of $\sqrt{-1}$ into the velocity field renders \eq{block-system} purely real.
The boundary conditions are
\Beq
\label{vr-bc} 
\left [\begin{array}{ccc} \left< 0,m+1,1\right| & \left< 0,m-1,1\right|  & 0 \end{array}\right] \ket{\Xi} \ = \ 0
\Eeq
This boundary-restriction operator replaces the last matrix row on top blocks of the left-hand side of \eq{block-system}; the blocks corresponding to $v^{+}$. A row of zeros replace the last matrix row of the top blocks on the right-hand side. We solve the numerical discretisation of \eq{block-system} with the same eigenvalue software as example~1. 

Like example~1, the inertial waves problem only allows one boundary condition for a system with (apparently) two radial derivatives. This results from the coordinate singularity causing a reduction in the differential order at $r=0$. This manifests in the matrix system because $D^{-}$ is a non-singular diagonal matrix, and hence requires no additional conditions to invert. Alternatively, the $D^{+}$ matrix contains a null space, and requires a boundary condition. We insert the radial-velocity boundary condition in the block corresponding to the $\ket{v_{m}^{+}}$ variable.  

Figure~\ref{eigenvalue-plot}b shows the relative eigenvalue error for the inertial wave problem with $m=1$, $\alpha=1$. Like the Bessel function example, there is a significant region with very high accuracy, and a region of the spectrum with significant error. Both of these regions occupy roughly the same proportions as the respective regions in the Bessel function problem. The major distinction is that the error occurs in the middle of the spectrum, rather than on one end. This results from the fact that small-horizontal-scale inertial wave oscillate slower than domain-filling modes. The inaccurate part of the spectrum results from under resolving small-scale modes in all cases.  In both example 1 \& 2, all the eigenvalues exist in the correct part of the complex plane; on the positive real axis for example 1, and on the real axis between $\pm1$ for example 2. While it is not possible to avoid inaccurate values for some modes, remaining on the real axis is not a guaranteed feature of a numerical approximation. Our scheme succeeds well in this respect. 

We discuss $m=1$ modes in this second example, but note that all modes preform well, just like in example 1. In the higher-$m$ cases the spectrum becomes squeezed to a smaller region between $\pm1$, otherwise there is little difference.  The eigenfunction error corresponds to the error in the eigenvalues. 

\paragraph{\textbf{Example 3}  (Linear perturbations of Hagen--Poiseuille pipe flow):} Hagen and Poiseuille first studied pipe flow experimentally in the 1830s \cite{sutera_skalak_1993}. This problem has guided studies investigating the onset of turbulence since the pioneering experiments of Reynolds in the 1880s \cite{reynolds_1883}. A lot of work continues today regarding the nonlinear dynamics of this system \cite{pringle_kerswell_2010}. After much investigation, a large sector of the turbulence-modelling community believes all linear perturbations in Poiseuille pipe flow eventually decay exponentially in time. Finite-amplitude perturbations likely produce the rich turbulent behaviour observed in laboratory experiments. This problem interests us because it poses a difficult numerical challenge; independent of the physical significance of the linear dynamics. Meseguer \& Trefethen (2003) \cite{meseguer2003linearized} discuss the linear problem at length. The eigenvalue problem for linear perturbations is:
\Beq
\label{vm-eq} \lambda  v_{m}(r)  + \nabla p_{m}(r) + L v_{m}(r) &=&  0 \n
\label{wm-eq} \lambda w_{m}(r) + W'(r) \,e_r\! \cdot v_{m}(r)   + i \alpha\, p_{m}(r) + L w_{m}(r) &=& 0  \n
\label{div-eq} \nabla \cdot v_{m}(r) + i \alpha\,  w_{m}(r) &=&  0.
\Eeq
where 
\Beq
\label{L-def} L \equiv i \alpha\, W(r) - \nu \nabla^{2} + \nu \alpha^2.
\Eeq 
For the background flow in the stream-wise direction, we choose the pressure-driven parabolic profile 
\Beq
W(r) = 1 - r^{2}, \quad \mathrm{with} \quad W'(r) = -2 r.
\Eeq
The field $v_{m}(r)$ represents the two-dimensional vector velocity in the pipe cross section. The scalar fields $w_{m}(r)$, and $p_{m}(r)$ represent the velocity in the stream-wise direction, and pressure respectively. 
The Laplacian implicit in \eq{L-def} represents the scalar Laplacian in \eq{wm-eq}, or the (more complicated) vector operator with curvature terms in \eq{vm-eq}. The advantage of the spinor approach is that these present the same level of complication in the numerical implementation. 
The parameters $\alpha,\lambda$ represent the mode wavenumber along the pipe axis, and the complex-valued growth rate respectively. The viscosity parameter $\nu = 1/\mbox{Re}$, where $\mbox{Re}$ represents the flow Reynolds number. We require no-slip boundary conditions on the outer wall
\Beq
v_{m}(r=1) = w_{m}(r=1) = 0.
\Eeq
The full eigenstate now contains four field components 
\Beq
\ket{\Xi_{m}} = \left[\begin{array}{cccc} \left| v_{m}^{+}\right> & \left|{v_{m}^{-}}\right> & \left|{w_{m}}\right> & \left|{p_{m}}\right> \end{array}\right]^{T}.
\Eeq
We solve the generalised eigenvalue problem for the complex-valued growth rate, $\lambda$,
\Beq
\mathcal{L} \ket{\Xi_{m}} \ = \ - \lambda\, \mathcal{R} \ket{\Xi_{m}}.
\Eeq
where 
\Beq
\label{LHS}
\mathcal{L} &\equiv&
\left[\begin{array}{cccc} L  & 0 & 0 & C D^{+} \\ 0 & L & 0 & C D^{-} \\ -\sqrt{2} CC R^{-} & -\sqrt{2} CC R^{+} &   L & i\alpha CC \\ D^{-} & D^{+}  & i\alpha C & 0  \end{array}\right],\n
\label{RHS} \mathcal{R} &\equiv&
\left[\begin{array}{cccc}  CC & 0 & 0 & 0 \\ 0 &  CC & 0 & 0 \\ 0 & 0 & CC &0 \\ 0 & 0 & 0 & 0  \end{array}\right],
\Eeq
and $L \equiv i \alpha CC (I-R^{+}R^{-})  - 2\nu  D^{+}D^{-} + \nu \alpha^{2} CC$.
The boundary conditions become
\Beq
\label{vr-vth-bc} 
\left [\begin{array}{cccc} \left< 0,m+1,1\right| & 0  & 0 & 0 \\ 
						0 & \left< 0,m-1,1\right| & 0 & 0 \\
						0 & 0 & \left< 0,m,1\right| & 0
 \end{array}\right] \ket{\Xi_{m}} \ = \ 0.
\Eeq
As with all other examples, we note that the operators in the different blocks of \eqs{LHS}{RHS} acquire different $(k,m)$ values depending on the field components. 

In addition to the eigenvalue problem, we solve for the parallel-component of the vorticity, $\omega_{3}$, and stream function, $\psi$. 
\Beq
\nabla^{2} \psi_{m} = \nabla \cdot \left( e_{3} \times v_{m}\right) = i ( \nabla_{-} v_{m}^{+} - \nabla_{+} v_{m}^{-}) \equiv \omega_{3},
\Eeq
which translates to
\Beq
2 D^{-}D^{+} \ket{\psi_{m}} = i C D^{-}\ket{v_{m}^{+}} -  i C D^{+}\ket{v_{m}^{-}}.
\Eeq
Often for an incompressible flow, the parallel velocity and stream function provide the dynamical variables in a higher-order-derivative scalar formulation of the original vector system. This example requires no special treatment to enforce the incompressibility constraint in \eq{div-eq}, which stands on equal footing with the other equations in the system.

\begin{figure}
\begin{center}
\includegraphics[scale=0.17]{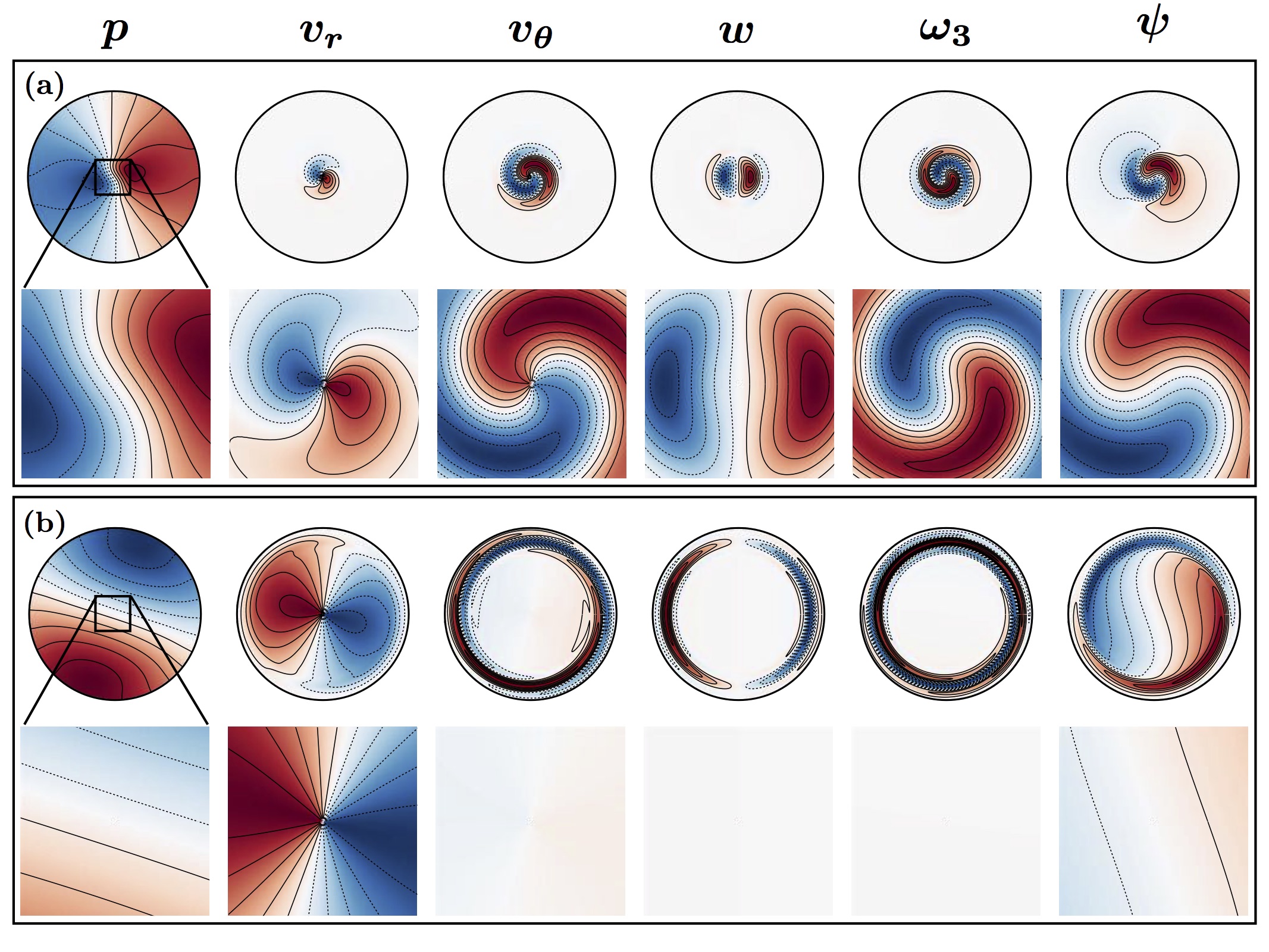}
\caption{Full-disk images for the first two modes for pipe-flow with $m=1$ and $\mathrm{Re}=3162$, and $N_{r}=200$ radial points. Panel (a) \& (b) shows the first- and second-slowest-decaying mode respectively. In each panel, the square plots show zoom-in plots for $|x|,|y| \le 0.2$. The behaviour at the origin remains regular for each field. 
\label{eigenmode-plot1}}
\end{center}
\end{figure}

\begin{figure}
\begin{center}
\includegraphics[scale=0.17]{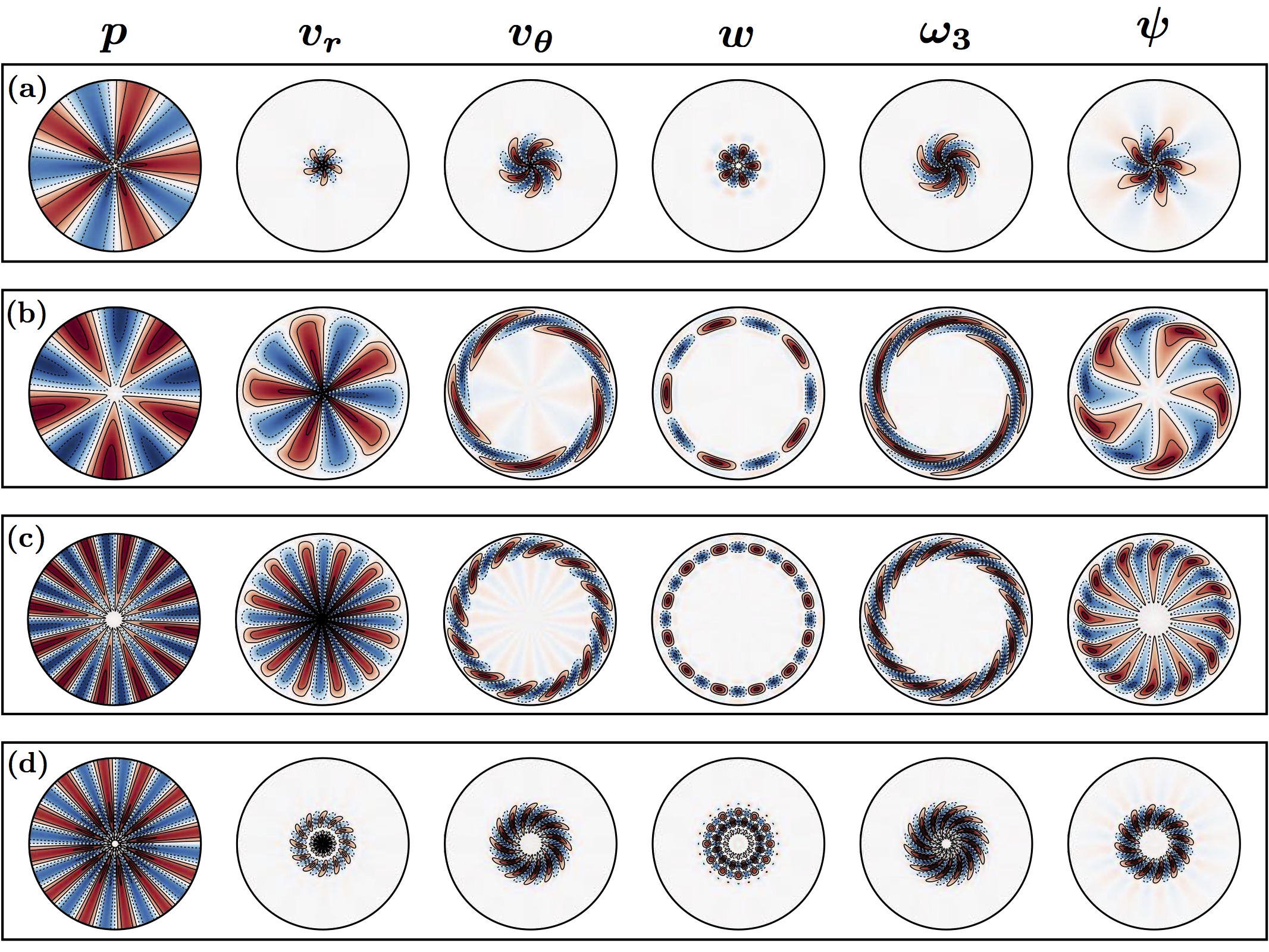}
\caption{Additional full-disk images for different pipe-flow modes with $\mathrm{Re}=3162$, and $N_{r}=200$ radial points. (a) corresponds to $m=5,n=1$, (b) corresponds to $m=5,n=2$, (c) corresponds to $m=12,n=1$, (d) corresponds to $m=12,n=5$. Zoom-in plot are omitted, but the behaviour at the origin remains smooth for each field.
\label{eigenmode-plot2}}
\end{center}
\end{figure}

Figure~\ref{eigenmode-plot1} shows full-disk eigenmodes for the first- and second- slowest decaying modes for $\mathrm{Re} = 10^{7/2}$, $\alpha=1$, and $m=1$. The plots show the four primitive dynamical variables $p,v_{r},v_{\theta},w$ along with the two \textit{ex post facto} derived fields $\omega_{3},\psi$. Below each full-disk image ($x^2+y^2 \le 1$), we show a zoomed-in square image on the region $\max(|x|,|y|) \le 0.2$. The zoom-in plots show no signs of numerical singularities for at the origin. The scalar fields, $p,w,\omega_{3}, \psi$ remain smooth, and the vector quantities show the correct parity. Figure~\ref{eigenmode-plot2} shows additional full-disk eigenmodes for $m=5$, and $m=12$. These plots show that the higher-$m$ still behave as required near the origin. Traditional methods (\eg Chebyshev bases) can cope with low-$m$ values, but fail for $m \gg 1$. Nonlinear simulations at high Reynolds number require a large range of azimuthal wavenumber to resolve the resulting turbulent states. For these situations, a scheme with a correct accounting for the origin is essential. 

Table~\ref{table:Pipe-Eigenvalues} shows a list of computed growth rates for different $m$ values for $\alpha=1$, and two different Reynolds numbers. In each case, we report the fastest growing (equivalently slowest decaying) modes of the centre- and wall-localised branches.
The $m=1$ values for both Reynolds numbers match the values reported in Meseguer \& Trefethen (2003) \cite{meseguer2003linearized} to the previously reported 11 digits. We tested all other values reported in \cite{meseguer2003linearized}. We match all cases. We report additional (higher $m$) values for potential future numerical comparison. The individual modes reported in table~\ref{table:Pipe-Eigenvalues} correspond to the same eigenmodes shown in figures~\ref{eigenmode-plot1}~\&~\ref{eigenmode-plot2}. For most cases, the first centre mode decays slower than the first wall mode. But this is not always the case. These modes can lie immediately adjacent to each other in the spectrum, but not necessarily so. For $\mathrm{Re}=10^4$, we converge to the reported solution after $\approx 50$ radial basis elements. For $\mathrm{Re}=10^7$ we converge after $\approx 200$ radial basis elements. In both cases the convergences is independent of $m$.

\begin{table}[ht]
\smallskip
\centering
\begin{tabular}{cccccc}
\hline
$m$ & type & $n$ & Re & Real($\lambda$) & Imag($\lambda$)  \\
\hline
\hline
1 & centre & 0 & $10^4$ & -0.0227049145535 & 0.951481194735  \\
\hline                      	        
1 & wall & 1& $10^4$ & -0.0472321995947 & 0.273788709331  \\
\hline                               
\hline                   
5 & centre & 0 & $10^4$ & -0.0725274157946 & 0.898561158159  \\
\hline
5 & wall & 1 & $10^4$ & -0.0793504734563 & 0.247410847332  \\
\hline
\hline
12 & wall & 0 & $10^4$ & -0.0948648867252  & 0.144951983763  \\
\hline
12 & centre & 3 & $10^4$ & -0.170456145014 & 0.800901547889  \\
\hline
\hline
1 & centre & 0 & $10^7$ & -0.000721091206991 & 0.998464685977  \\
\hline                                   
1 & wall & 15 & $10^7$ & -0.00748956875998 & 0.0303389812102  \\ 
\hline   
\hline                     
5 & centre & 0 & $10^7$ & -0.00229096203822 & 0.996790918537  \\
\hline 
5 & wall & 15 & $10^7$ & -0.00855398926555 & 0.0148836399355  \\
\hline  
\hline                         
12 & centre & 0 & $10^7$ & -0.00538731680888 & 0.993703412087  \\
\hline                           
12 & wall & 5 & $10^7$ & -0.00784725003139 & 0.0296167267785  \\
\hline
\hline
\end{tabular}
\caption{Pipe-flow eigenvalues computed in example 3. In each case, the ``centre'' and ``wall'' modes represent the slowest-decaying solution with eigenfunctions localised near the origin or near the boundary, respectively. The number $n$ represents where a given mode occurs in a sorted list of growth rates.}
\label{table:Pipe-Eigenvalues}
\end{table}

\begin{figure}
\begin{center}
\includegraphics[scale=0.35]{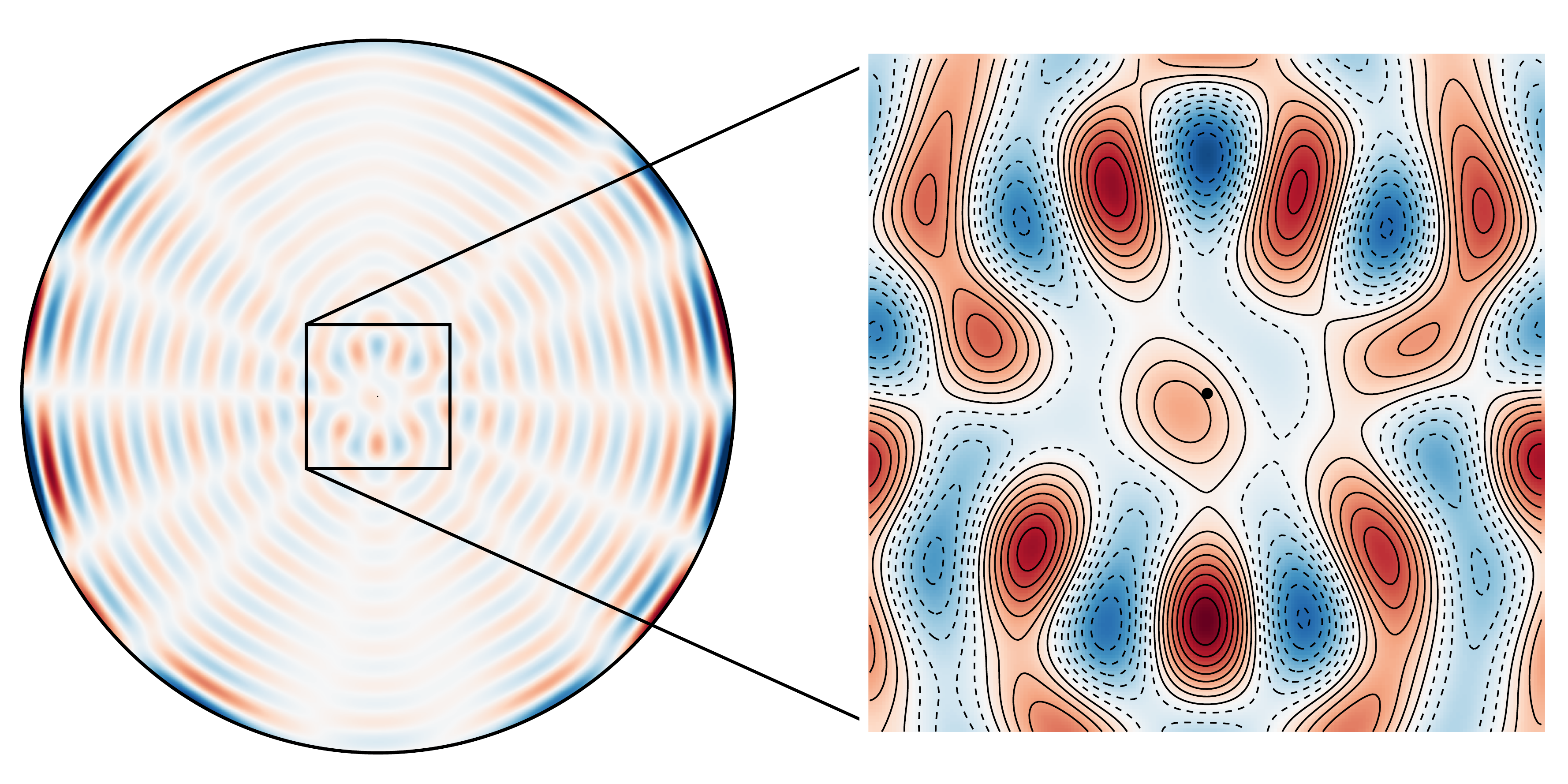}
\caption{Contour plot of the calculated solution to the forced Helmholtz equation for $\kappa=60$, with $s(x,y)=e^{-(x-0.4)^2-(y-0.3)^2}$ and $g(x,y)=y\cos(10x)$.  No artificial singularity appears near the origin (marked with small dot in the zoomed-in image). 
\label{helmholtz-plot}}
\end{center}
\end{figure}

\paragraph{\textbf{Example 4} (Forced Helmholtz equation):}  We consider the forced Helmholtz equation 
\Beq
\label{Helmholtz}
\left(\nabla^{2}  + \kappa^2\right) f = s \qquad\hbox{and}\qquad  f(r=1,\theta) = g(\theta)
\Eeq
with $s(x,y)=e^{-(x-0.4)^2-(y-0.3)^2}$ and $g(x,y)=y\cos(10x)$.  When  $\kappa$ is very  large, the solutions become increasingly oscillatory; see Figure~\ref{helmholtz-plot}. The zoom-in plot on the right demonstrates that no artificial singularity arises at the origin, despite the complicated nature of the solution. We choose this problem to demonstrate the performance of the methods on a challenging and natural problem requiring simultaneous $r,\theta$ dependence. 

\begin{figure}
\begin{center}
\includegraphics[scale=0.30]{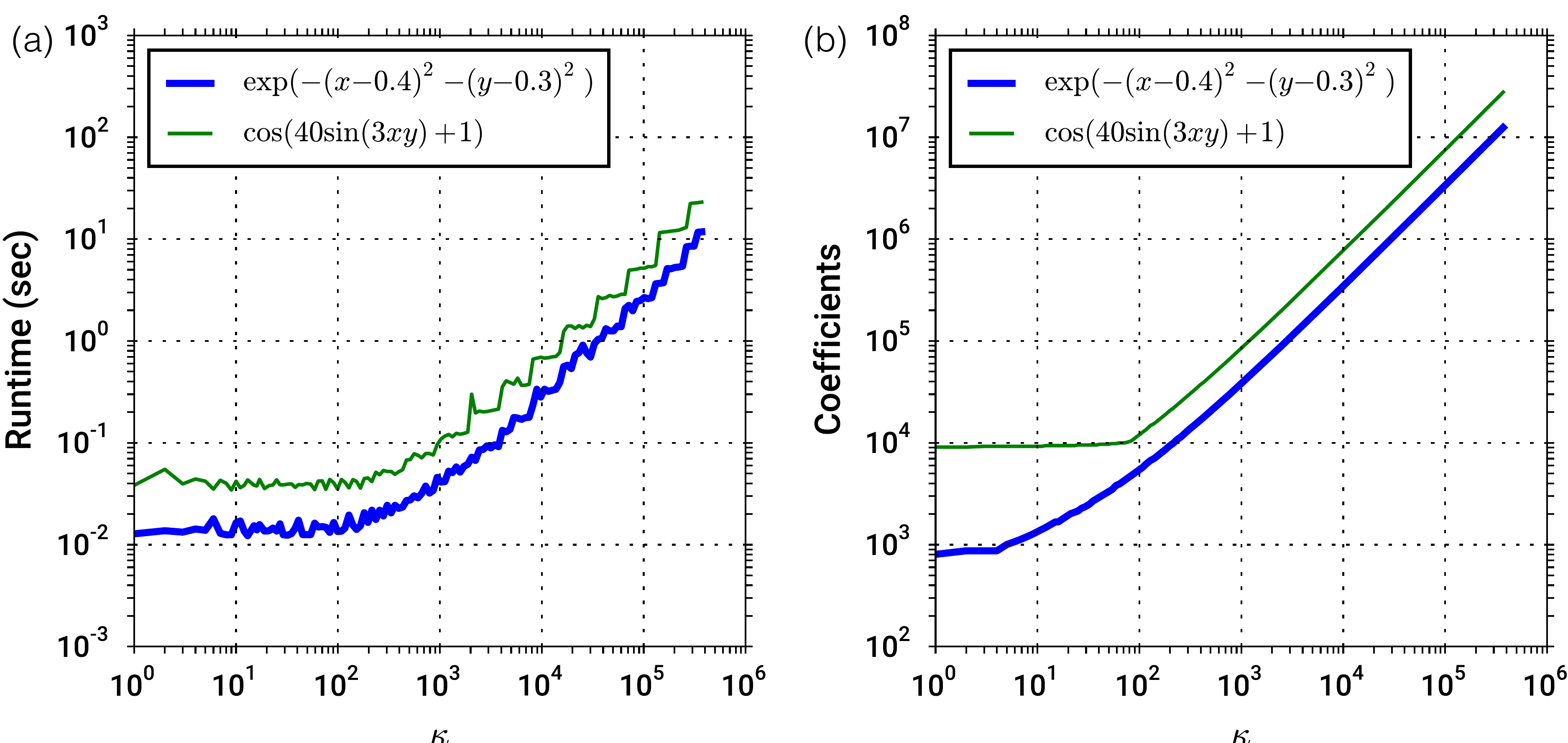}
\caption{(a) Time taken to calculate solutions to the Helmholtz equation for two different forcing terms (green versus blue) with $g(x,y)=y\cos(10x)$.  (b) Total number of coefficients  needed to resolve  the solutions. The asymptotic slope becomes linear in $\kappa$ after an initial fixed costs needed to resolve the right-hand side forcing, as indicated by the background dashed grid. 
\label{helmholtz-time}}
\end{center}
\end{figure}

The bandedness of our representation implies the complexity of the solution is $O(M N )$ operations, where $M$ is the number of $\theta$ modes needed to resolve $g(\theta)$ and $s(r,\theta)$ in the $\theta$ direction, and $N$ is the optimal number of polynomial coefficients needed to solve to a given tolerance in radius.  In Figure~\ref{helmholtz-time} we demonstrate the linear complexity for two different resolved right-hand sides. We approximate, $e^{-(x-0.4)^2-(y-0.3)^2}$ with 25 $\theta$ modes and 11 radial polynomial coefficients per mode, while $\cos(40 \sin(3xy)+1)$ requires 149 $\theta$ modes and  60 polynomial coefficients per mode.   We use the adaptive QR method implemented in the ApproxFun.jl Julia package \cite{olver2013fast,olver2014practical}. This adaptively determines the optimal number, $N$, of polynomial coefficients needed to resolve the solution for each $\theta$ mode.  \Eq{Helmholtz} presents a challenge because the solution can become much more structured than the forcing for $\kappa \gg 1$. At the high rage of $\kappa$, we solve an equation with over 25 million unknowns in less than 25~seconds\footnote{Timings on a 2012 MacBook Pro, with an 2.7 GHz Intel Core i7.}. 

We also solve a screened Poisson equation similar in form the \eq{Helmholtz} only with $\kappa^{2} < 0$. Rather than small-scale oscillations, the solutions tend to form sharp boundary layers near $r=1$. The performance efficiency and accuracy of the methods treats this case equally as well as the Helmholtz problem. 

\section{Conclusion}

This paper produces a new set of algorithms for numerical computations on the unit disk. Choosing a hierarchy of bases built out of Jacobi polynomials allows for solving a range of important differential equations. We show a range of time-independent examples. Extending to time-evolving system only requires  conducting the type of calculations presented here in the iterative fashion. 

The methods described here explicitly construct basis elements that respect the coordinate singularity in polar coordinates at $r=0$. Beyond the disk, the type of calculus presented here will exist for other geometries with coordinate singularities. The two-dimensional disk represents only a prototypical example.  

The Laguerre class of polynomials allows representing functions in polar coordinates for the whole two-dimensional plane. In particular, 
\Beq
\lim_{R \rightarrow \infty}r^{m} P_n^{(R^2,m)}(2(r/R)^2-1)(1-(r/R)^2)^{R^{2}}  = (-1)^n r^{m} L_n^{(m)}(r^2)e^{-r^2}.
\Eeq
Many of the formulae present in this paper extend to the infinite plane for $k\to \infty$ with only minor adjustments for the Gaussian factor. In this limit, the conversion matrix becomes the identity, $C \to I$.  

Two-dimensional spheres and solid balls also display similar behaviour as the unit disk. The sphere effectively contains a local disk at both poles. Spherical harmonic functions naturally resolve the coordinate singularities at the poles for scalar functions. However, a wider class of Jacobi polynomials allows the use of sparse vector calculus for higher-order tensor fields. In addition to the poles, the centre of a solid sphere requires relations completely analogous to \eqs{analytic}{v-analytic} only for the spherical harmonic degree \cite{orszag_1974}. In an upcoming publication, we show how to organise sparse tensor calculus in three-dimensional spheres. This includes algebraic structures similar to those found here for the disk. 

The algebraic-analytic link is an important theme this current paper. As a general principle, the singularities of a PDE, either coordinate-based or physical, dictate what types of bases are needed to represent meaningful solutions. When followed carefully, this guiding principle should allow the construction of numerical methods more efficient and straightforward to implement than more traditional schemes. Spherical geometry exhausts the parameter variation found in Jacobi polynomials. However, many more lesser-explored sets of orthogonal polynomials exist, with a large variety of symmetry properties. This includes discrete sets. We encourage the future development of numerical methods based on the guiding algebraic-geometric properties of a given general problem class. 

\textbf{ACKNOWLEDGMENTS}

GMV acknowledges support from the Australian Research Council, project number DE140101960. KJB is supported by a National Science Foundation Graduate Research Fellowship under Grant No.~1122374. DL is supported by the Hertz Foundation.  SO acknowledges support from the Australian Research Council, project number DE130100333

\appendix
\renewcommand*{\thesection}{\Alph{section}}

\section{Matrix entries}

In this appendix $\delta_{i,j}$ represents the Kronecker delta. We show all matrix coefficients as they act on columns of spectral coefficients from the left. 

\paragraph{Differentiation matrices}

\Beq
D^{+}: \{k,m\} \rightarrow \{k+1,m+1\}
\Eeq
\Beq
D^{+}_{n',n} = \sqrt{2n(n+k+m+1)} \, \delta_{n', n-1}	
\Eeq

\Beq
D^{-}: \{k,m\} \rightarrow \{k+1,m-1\}
\Eeq
\Beq
D^{-}_{n',n} = \sqrt{2(n+k+1)(n+m)} \, \delta_{n', n}	
\Eeq

\paragraph{Multiplication matrices}

\Beq
R^{+}: \{k,m\} \rightarrow \{k,m+1\}
\Eeq
\begin{align}
\label{R+ elements}
R^{+}_{n',n} = \sqrt{\frac{(n+m+1)(n+k+m+1)}{(2n+k+m+1)(2n+k+m+2)}} \, &\delta_{n',n} \\ + \sqrt{\frac{n(n+k)}{(2n+k+m)(2n+k+m+1)}} \, &\delta_{n',n-1}
\end{align}

\Beq
R^{-}: \{k,m\} \rightarrow \{k,m-1\}
\Eeq
\begin{align}
\label{R- elements}
R^{-}_{n',n} = \sqrt{\frac{(n+1)(n+k+1)}{(2n+k+m+1)(2n+k+m+2)}} \, &\delta_{n',n+1} \\ + \sqrt{\frac{(n+m)(n+k+m)}{(2n+k+m)(2n+k+m+1)}} \, &\delta_{n',n}
\end{align}

\paragraph{Conversion matrices}

\Beq
C: \{k,m\} \rightarrow \{k+1,m\}
\Eeq
\begin{align}
\label{C elements}
C_{n',n} = \sqrt{\frac{(n+k+1)(n+k+m+1)}{(2n+k+m+1)(2n+k+m+2)}} \, &\delta_{n',n} \\ - \sqrt{\frac{n(n+m)}{(2n+k+m)(2n+k+m+1)}} \, &\delta_{n',n-1}	
\end{align}

\begin{equation}
C^{\dagger}: \{k,m\} \rightarrow \{k-1,m\}
\end{equation}
\begin{align}
\label{Cdag elements}
C^{\dagger}_{n',n} = \sqrt{\frac{(n+1)(n+m+1)}{(2n+k+m+1)(2n+k+m+2)}} \, &\delta_{n',n+1} \\ + \sqrt{\frac{(n+k)(n+k+m)}{(2n+k+m)(2n+k+m+1)}} \, &\delta_{n',n}
\end{align}

\section{Dirac's bra-ket notation}

Dirac notation is uncommon in numerical analysis applications outside of quantum mechanics.  We believe the notation provides many advantages with regard to clarity and abstraction. This short appendix reviews the important aspects of the notation as it applies to the bulk of the paper. 

In the appendix, $Q_{n}(r)$ denotes any one of the Jacobi-based orthogonal polynomials in the main section. We drop the $k,m$ indices that parameterise different bases.
We represent functions of $r$ such that  
\Beq
f(r) \ = \ \sum_{n=0}^{\infty} f_{n} Q_{n}(r).
\Eeq
For a given weight function $w(r)$ [\eg $r(1-r^2)^k$],
\Beq
\int_0^1 Q_{n}(r) Q_{n'}(r) w(r) \mathrm{d} r = \delta_{n,n'}.
\Eeq
Therefore a function's spectral-expansion coefficients result from projecting against a given polynomial  
\Beq
 f_{n} \ = \ \int_0^1 f(r) Q_{n}(r)  w(r) \mathrm{d} r.
\Eeq

We define a (discrete) complete orthonormal basis such that
\Beq
\braket{n}{n'} \ = \ \delta_{n,n'}, \quad 
\sum_{n=0}^{\infty} \ket{n}\bra{n} \ = \ \mathrm{I}.
\Eeq
The element $\bra{n}$ denotes an infinite row vector with $0$'s in the first $n$ entries, a $1$ in the $(n+1)$st entry, and $0$'s after that. For example, 
\Beq
\bra{0} = [\,1, 0, 0, 0, 0, 0, \ldots\,],\quad \bra{4} = [\,0, 0, 0, 0, 1, 0, \ldots\,]. 
\Eeq

In terms of the discrete basis and orthonormal functions, we define the continuous basis element corresponding to evaluation at the point $r$
\Beq
\bra{r}  \ \equiv  \ \sum_{n=0}^{\infty} Q_{n}(r) \bra{n} 
\Eeq
\Beq
\int_{0}^{1} \ket{r} \bra{r}  w(r) \mathrm{d} r \ = \ \sum_{n=0}^{\infty}\sum_{n'=0}^{\infty}\ket{n} \bra{n'} \int_{0}^{1}Q_{n'}(r) Q_{n'}(r) w(r) \mathrm{d} r 
\Eeq
Therefore,
\Beq
\int_{0}^{1} \ket{r} \bra{r} w(r) \mathrm{d} r \ = \  \mathrm{I}.
\Eeq
Also, 
\Beq
\braket{r}{r'}  \ \equiv  \ \sum_{n=0}^{\infty} Q_{n}(r) Q_{n}(r')  \ = \ \frac{1}{w(r)}\delta(r-r'), 
\Eeq
which follows from the Christoffel--Darboux formula for orthogonal polynomials.

The orthogonal functions represent the coefficients of the $r$ basis in terms of the $n$ basis, 
\Beq
\braket{r}{n} = Q_{n}(r).
\Eeq
Now,
\Beq
f(r) = \braket{r}{f} = \sum_{n=0}^{\infty} Q_{n}(r) \braket{n}{f} 
\Eeq
where
\Beq
f_{n} = \braket{n}{f}
\Eeq
represent the coefficients of the state $f$ in the $n$ basis. 

The $L^{2}$ inner product of functions takes the form 
\Beq
\int_{0}^{1} g^{*}(r)f(r) w(r) \mathrm{d} r = \int_{0}^{1} \braket{g}{r} \braket{r}{f} w(r) \mathrm{d} r  
\Eeq
therefore 
\Beq
\int_{0}^{1} g^{*}(r)f(r) w(r) \mathrm{d} r =  \braket{g}{f}.
\Eeq


\section*{References}

\begin{thebibliography}{10}

\bibitem{bhatia1954circle}
{\sc A.~Bhatia and E.~Wolf}, {\em On the circle polynomials of Zernike and
  related orthogonal sets}, in Mathematical Proceedings of the Cambridge
  Philosophical Society, vol.~50, Cambridge University Press, 1954, pp.~40--48.

\bibitem{boyd_book}
{\sc J.~P. Boyd}, {\em Chebyshev and Fourier Spectral Methods: Second Revised Edition}, Dover (2001)

\bibitem{boyd2011comparing}
{\sc J.~P. Boyd and F.~Yu}, {\em Comparing seven spectral methods for
  interpolation and for solving the poisson equation in a disk: Zernike
  polynomials, Logan--Shepp ridge polynomials, Chebyshev--Fourier series,
  cylindrical Robert functions, Bessel--Fourier expansions, square-to-disk
  conformal mapping and radial basis functions}, Journal of Computational
  Physics, 230 (2011), pp.~1408--1438.

\bibitem{coutsias_hagstrom_1996}
{\sc E.~A.~Coutsias, T.~Hagstrom and D.~Torres}, {\em An efficient spectral method for ordinary
differential equations with rational function coefficients}, Mathematics of Computation, 65 (1996), pp. 611--635.

\bibitem{doha2002efficient}
{\sc E.~H. Doha and W.~M. Abd-Elhameed}, {\em Efficient spectral-Galerkin
  algorithms for direct solution of second-order equations using ultraspherical
  polynomials}, SIAM Journal on Scientific Computing, 24 (2002), pp.~548--571.

\bibitem{doha2006efficient}
{\sc E.~H. Doha and A.~H.~Bhrawy}, {\em Efficient spectral-Galerkin algorithms for direct solution for second-order differential equations using Jacobi polynomials}, Journal of Numerical Algorithms, 42 (2006),  pp.~137--164

\bibitem{doha2009efficient}
{\sc E.~Doha and W.~Abd-Elhameed}, {\em Efficient spectral
  ultraspherical-dual-Petrov--Galerkin algorithms for the direct solution of
  $(2n+ 1)$th-order linear differential equations}, Mathematics and Computers in
  Simulation, 79 (2009), pp.~3221--3242.

\bibitem{dunklxu}
{\sc C.~F. Dunkl and Y.~Xu}, {\em Orthogonal polynomials of several variables,
  Second Edition}, Cambridge University Press, 2014.

\bibitem{dirac_1939}
{\sc P.A.M.~Dirac}, {\em  A new notation for quantum mechanics}, Mathematical Proceedings of the Cambridge Philosophical Society, 35 (1939), pp.~1--4.

\bibitem{fornberg1998practical}
{\sc B.~Fornberg}, {\em A practical guide to pseudospectral methods}, vol.~1,
  Cambridge University Press, 1998.

\bibitem{gardner_etal_1989}
{\sc D.~R.~Gardner, S.~A. Trogden, and R.~W.~Douglass}, {\em A modified tau spectral method that eliminates spurious eigenvalues}, Journal of Computational Physics, 80 (1989), pp.~137--167?

\bibitem{greengard_1991}
{\sc L.~Greengard}
{\em Spectral integration and two-point boundary value problems} SIAM Journal of Numerical Analysis, 28 (1991), pp.~1071--1080

\bibitem{greenspan_1968}
{\sc H.~P.~Greenspan}, {\em The theory of rotating fluids}, vol.~1,
  Cambridge University Press, 1968.

\bibitem{julien_watson_2009}
{\sc K.~Julien, M.~Watson}, {\em Efficient multi-dimensional solution of PDEs using Chebyshev spectral methods}, Journal of Computational Physics, 228 (2009), pp.~1480--1503.

\bibitem{koornwinder}
{\sc T.~Koornwinder}, {\em Two-variable analogues of the classical orthogonal
  polynomials}, Theory and applications of special functions,  (1975),
  pp.~435--495.

\bibitem{li2010optimal}
{\sc K.~Li, P.~W. Livermore, and A.~Jackson}, {\em An optimal Galerkin scheme
  to solve the kinematic dynamo eigenvalue problem in a full sphere}, Journal
  of Computational Physics, 229 (2010), pp.~8666--8683.

\bibitem{livermore_jones_worland_2007}
{\sc P.~W. Livermore, C.~A.~Jones and S.~J.~Worland}, {\em Spectral radial basis functions for full sphere computations}, Journal of Computational Physics, 227 (2007), pp. 1209--1224

\bibitem{matsushima1995spectral}
{\sc T.~Matsushima and P.~Marcus}, {\em A spectral method for polar
  coordinates}, Journal of Computational Physics, 120 (1995), pp.~365--374.

\bibitem{mcfadden_etal_1990}
{\sc G.~B. McFadden, B~T. Murray, and R~F.~Boisvert}, {\em Elimination of spurious eigenvalues in the Chebyshev tau spectral method}, Journal
  of Computational Physics, 91 (1990), pp.~228--239.

\bibitem{meseguer2003linearized}
{\sc A.~Meseguer and L.~N. Trefethen}, {\em Linearized pipe flow to Reynolds
  number $10^7$}, Journal of Computational Physics, 186 (2003), pp.~178--197.

\bibitem{muite_2010}
{\sc B.~K.~Muite}
{\em A numerical comparison of Chebyshev methods for solving fourth order semilinear initial boundary value problems}, Journal of Computational and Applied Mathematics, 234 (2010), pp. 317--342

\bibitem{DLMF}
{\sc F.~W.~J. Olver, D.~W. Lozier, R.~F. Boisvert, and C.~W. Clark}, {\em {NIST
  Handbook of Mathematical Functions}}, Cambridge University Press, 2010.

\bibitem{olver2013fast}
{\sc S.~Olver and A.~Townsend}, {\em A fast and well-conditioned spectral
  method}, SIAM Review, 55 (2013), pp.~462--489.

\bibitem{olver2014practical}
{\sc S.~Olver and A.~Townsend}, {\em A practical framework for
  infinite-dimensional linear algebra}, in HPTCDL, 2014, pp.~57--62.

\bibitem{orszag_1971}
{\sc S.~A.~Orszag}, {\em Accurate solution of the Orr-Sommerfeld stability equation}, Journal of Fluid Mechanics, 50 (1971), pp.~689--703.

\bibitem{orszag_1974}
{\sc S.A.~Orszag}, {\em Fourier series on spheres}, Monthly Weather Review, 102 (1974), pp.~56--75.

\bibitem{pringle_kerswell_2010}
{\sc C.~C.~T.~Pringle and R.~R.~Kerswell}, {\em Using nonlinear transient growth to construct the minimal seed for shear flow turbulence}, Physical Review Letters, 105 (2010), 154502

\bibitem{reynolds_1883}
{\sc O.~Reynolds}
{\em An experimental investigation of the circumstances which determine whether the motion of water shall be direct or sinuous, and of the law of resistance in parallel channels}, Philosophical Transactions of the Royal Society of London, 174 (1883), pp. 935--982 

\bibitem{sutera_skalak_1993}
{\sc S.~P.~Sutera, and R.~Skalak}, {\em The history of Poiseuille's law},
  Annual Review of Fluid Mechanics, 25 (1993), pp.~1--19.

\bibitem{slevinsky_2016}
{\sc R.~M.~Slevinsky}, {\em On the use of Hahn's asymptotic formula and stabilized recurrence for a fast, simple, and stable Chebyshev--Jacobi transform} (2016) \url{arxiv.org/pdf/1602.02618.pdf}

\bibitem{sakai2009application}
{\sc T.~Sakai and L.~Redekopp}, {\em An application of one-sided Jacobi
  polynomials for spectral modelling of vector fields in polar coordinates},
  Journal of Computational Physics, 228 (2009), pp.~7069--7085.

\bibitem{townsend_olver_2015}
{\sc A.~Townsend and S.~Olver}, {\em The automatic solution of partial differential equations using a global spectral method}. Journal of Computational Physics, 299 (2015), pp. 106--123.

\bibitem{townsend_etal_2016}
{\sc A.~Townsend, H.~Wilber, and G.~Wright}, {\em Computing with functions in spherical and polar geometries II. The disk}. (2016) \url{http://arxiv.org/abs/1604.03061}.

\bibitem{trefethenspectralmethods}
{\sc L.~N. Trefethen}, {\em Spectral methods in MATLAB}, vol.~10, Siam, 2000.

\bibitem{viswanath_2015}
{\sc D.~Viswanath}, {\em Spectral integration of linear boundary value problems}, Journal of Computational and Applied Mathematics, 290 (2015), pp.~159--173


\bibitem{zernike}
{\sc F.~Zernike}, {\em Beugungstheorie des schneidenverfahrens und seiner
  verbesserten form, der phasenkontrast-methode}, Physica, 1,  (1934),
  pp.~689--704.

\bibitem{zebib_1984}
{\sc A.~Zebib}, {\em A Chebyshev method for the solution of boundary value problems}, Journal of Computational Physics,  53 (1984), pp.~443--455.

\end{thebibliography}
\end{document}